\documentclass[12pt]{article}
\def\Zp{${\bZ}_p$}  
 \def\row#1{(#1_0,\ldots,#1_{p-1})}
 
\def\xhat{{\hat x}} 
\def\k{\omega}
\def\Cp{${\bC}^p$}  

\def\0p{$({\bC}\setminus\{0\})^p\ $}  
\def\0n{$({\bC}\setminus\{0\})^n\ $}  
\def\C0pmat{({\bC}\setminus\{0\})^p\ }  
\def\C0nmat{({\bC}\setminus\{0\})^n\ }  
\def\eqref#1{{\rm(\ref{#1})}}
\def\tfrac#1#2#3{{\textstyle{#1\over #2}}#3}

\usepackage{a4wide,latexsym,amssymb}

\newtheorem{theorem}{Theorem}[section]
\newtheorem{lemma}[theorem]{Lemma}
\newtheorem{proposition}[theorem]{Proposition}
\newcommand{\rank}{\mathrm{rank}}

\newcommand{\Ker}{\mathrm{Ker}}
\newcommand{\Span}{\mathrm{span}}

\newcommand{\Arccos}{\mathrm{Arccos}}
\newcommand{\Arccot}{\mathrm{Arccot}}

\newcommand{\sign}{\mathrm{sign}}

\newcommand{\range}{\mathrm{range}}

\newtheorem{remark}[theorem]{Remark}
\newcommand{\bproof}{\noindent{\bf Proof: }}

\newcommand{\bC}{{\mathbb C}}
\newcommand{\bZ}{{\mathbb Z}}
\newcommand{\bN}{{\mathbb N}}
\newcommand{\bR}{{\mathbb R}}

\renewcommand{\theequation}{\arabic{section}.\arabic{equation}}
\def\varepsilon{\epsilon}

\begin{document}
\renewcommand{\theequation}{\arabic{section}.\arabic{equation}}
\title{All cyclic $p$-roots of index 3, found by symmetry-preserving
  calculations} \author{G\"oran Bj\"orck\footnote{Department of 
  Mathematics, Stockholm University, SE-106 91 $\,$ STOCKHOLM, 
  Sweden,  bjorck@math.su.se}$\,\,$ and 
  Uffe Haagerup\footnote{Department of Mathematics and Computer 
  Science, University of Southern Denmark, Campusvej 55, DK-5230 
  Odense M, Denmark,  haagerup@imada.sdu.dk}} \maketitle
\setcounter{section}{0}
\section*{Introduction}
\label{sec-introduction}
When using a Groebner basis to solve the highly symmetric system of
algebraic equations defining the cyclic $p$-roots, one has the feeling
that much of the advantage of computerized symbolic algebra over hand
calculation is lost through the fact that the symmetry is immediately
``thrown out'' by the calculations. In this paper, the problem of
finding (for all relevant primes $p$) all cyclic $p$-roots of index 3
(as defined in Section \ref{sec-notation}) is treated with the
symmetry preserved through the calculations. Once we had found the
relevant formulas, using MAPLE and MATHEMATICA, the calculations
could even be made by hand. On the other hand, with respect to 
a straightforward attack with Groebner basis, it is not even clear
how this could be organized for a general $p$.


In other terminologies, our results involve listings of all 
bi-unimodular sequences constant on the cosets of the group $G_{0}$
of cubic residues, or equivalently all circulant complex Hadamard 
matrices related to $G_{0}$ (cf.~\cite{CR}).

The corresponding problem for bi-unimodular sequences of index 2 was 
solved by the first named author in \cite{pisa} and shortly after solved 
independently by de la Harpe and Jones \cite{jones} in the case 
$p\equiv 1\pmod 4$  and by
Munemasa and Watatani \cite{munemasa} in the case $p\equiv 3\pmod{4}$  , see also 
\cite{haagerup}, sect. 3.

The organization of the paper should be clear from the section 
headings with the understanding that ``the main problem'' refers to {\it simple} 
sequences of index 3 (cf. Definitions \ref{simple}, 
 \ref{general}, and \ref{cycind3}).

\section{Notation, definitions, and problem formulation}
\label{sec-notation}
\setcounter{equation}{0}
We begin by quoting from ~\cite{pisa} and ~\cite{CR} definitions of
and relations between {\it bi-unimodular $p$-sequences} and {\it
  cyclic $p$-roots} for any positive integer $p$.
For any $p$-sequence $x$, that is any sequence $x=\row x$ of
$p$ complex numbers, define its normalized Fourier
transform by 
$\xhat_{\nu} ={1\over \sqrt{p}} \sum_{j=0}^{p-1}x_j\k^{j\nu}$, where
$\k=\exp ({2\pi i\over p}).$ 
The sequence $x$ is called {\it unimodular} if $|x_j|=1$ for
$j=0,1,\ldots,p-1$, ~and it is called {\it bi-unimodular} if both $x$
and $\xhat$ are unimodular.

Taking all indices modulo $p$, we define the {\it periodic
  autocorrelation coefficients} $\gamma_k$ by
\begin{equation} \gamma_k=\sum_{j\pmod p} {\bar x}_jx_{j+k}.
\label{eq:corr}
\end{equation}
Then, by the Parseval relation and an easy calculation,
\begin{equation} \xhat~ {\rm is~ unimodular}~ \Leftrightarrow
  (\gamma_0=p~ {\rm and}~ \gamma_k=0~{\rm when}~ k\not \equiv 0 \pmod
  p).\label{eq:hatuni}
\end{equation}
We will now express the property of bi-unimodularity with the help of
a certain system of algebraic equations.  Let $z=\row z \in {\bC}^p$. We will call 
$z$ a {\it ``cyclic $p$-root''}, if $z$
satisfies the following system of $p$ algebraic equations:
\begin{eqnarray}
  z_0+z_1+\cdots +z_{p-1} & = & 0, \nonumber \\ z_0z_1+z_1z_2+\cdots
  +z_{p-1}z_0 & = & 0, \nonumber\\ & \vdots & \label{eq:cyk}\\ 
  z_0z_1\cdots z_{p-2}+z_1z_2\cdots z_{p-1}+\cdots +z_{p-1}z_0\cdots
  z_{p-3} & = & 0, \nonumber\\ z_0z_1\cdots z_{p-1} & = & 1.\nonumber
\end{eqnarray}
(Note that the sums are cyclic and contain just $p$ terms and are in
general {\it not} the elementary symmetric functions.)  Let now
$x\in{\bC}^p$ and $z\in{\bC}^p$ be related by\begin{equation}
z_j=x_{j+1}/x_j
\label{eq:xtoz}\end{equation}
(with $x_p:=x_0$). Clearly $x$ is unimodular iff $ {\bar x}_j=1/x_j
(\forall j)$. In this case, (\ref{eq:corr}) for $k=1,2,\ldots,p-1$ 
becomes the $k$'th equation of (\ref{eq:cyk}).  Let us call $x$ 
{\it normalized} if
 $x_0=1$. Then \eqref{eq:hatuni} can be expressed as follows:
\newtheorem  {xiffz}{Proposition}[section]
\begin{xiffz}
A normalized $x=(1,x_{1},x_{2},\ldots,x_{p-1})$ is 
  bi-unimodular if and only if the corresponding $z$ is a {\rm
  unimodular} cyclic $p-$root.
\label{xiffz}
\end {xiffz}
In the rest of the paper, $p$ will be a {\it prime} $\equiv 1
\pmod{6}$, and we will define $s:=(p-1)/3$. The multiplicative group  
${\bZ}^*_p$ on ${\bZ}_p\setminus \{0\}$ is cyclic (cf. \cite{CR}) and 
has a unique index-3 subgroup $G_{0}$ (the group of cubic residues 
modulo $p$). Let $G_1$ and $G_2$ be the other two cosets of $G_0$ in 
$\bZ^*_p$. (The choice of the subscripts 1 and 2 will be
specified later.). We will now for $p$-sequences define a property 
{\it  of index} 3 meaning ``taking few values in a way governed 
by $G_0$'':

\newtheorem  {simple}[xiffz]{Definition}
\begin{simple}
\label{simple}
We will say that $x\in{\bC}^p$ is {\rm simple of
  index 3}, if there are complex numbers, $c_0,c_1,{\rm~
  and}~c_2$, such that
\begin{equation} 
  x_j=c_k~ {\rm when}~0\ne j\in G_k~~~(k=0,1,2).
\label{ind3}\end{equation}
\end{simple}

Note that we have slightly changed the notation from \cite{pisa} 
where index 3 was called ``pre-index 3'' and where ``index 3'' 
excluded the case of index 1. i.e. $c_{0}=c_{1}=c_{2}.$

Allowing shifts and multiplication by exponentials in a way familiar
in Fourier transform theory, we make the following definition:
\newtheorem  {general}[xiffz]{Definition}
\begin{general}
\label{general}
We will say that $x\in$\Cp~ has {\rm
  index 3}, if for some fixed elements $r\ne0$ and $l$ of $\bZ_p$ and some 
  {\rm
  simple} $y$ of index {\rm 3} we have
\begin{equation}
  x_j=\k^{rj}y_{j-l},
\end{equation}
which amounts to
\begin{equation}
  x_j=\k^{rj}c_k~ {\rm when}~0\ne j-l\in G_k~~~(k=0,1,2).
\label{genx}
\end{equation}
\end{general}

We will now define simple and general {\it cyclic $p$-roots} of index 3:
\newtheorem  {cycind3}[xiffz]{Definition}
\begin{cycind3}
\label{cycind3}
By a cyclic $p$-root of index {\rm 3} we will
mean a cyclic $p$-root $z$ such that the corresponding $x$, as defined
by {\rm(\ref{eq:xtoz})} has index {\rm 3}. We will also call a cyclic $p$-root $z$ 
 simple of index {\rm 3}, if the corresponding $x$ is simple of index {\rm 3}. 
\end{cycind3}
Note that we do not require $x$ (and thus $z$) to be unimodular.

The purpose of the present paper is to find explicitly all cyclic $p$-roots of index 3
(for every relevant prime $p$) using a method which utilizes the
symmetries of the system.

We will now show (following {\cite{pisa}), that if 
  $z$ is a simple cyclic $p$-root of
  index 3 and its corresponding $x$ is normalized by $x_0=1$,
  then the system (\ref{eq:cyk}) 
  reduces to a system of three equations for $c_0, c_1~{\rm and}~
  c_2$. (To help the reader, an example is given at the end of the
  section.) Let $g$ be a generator for $\bZ^*_p$, and let
  $G_0,G_1,G_2$ be the cosets of $G_0$, numbered in such a way that
  $G_k=\{g^{k+3m};m=0,1,\ldots,s-1\}.$ For every $i$
  and $k=0,1,2$, and every $d=1,\ldots,p-1$, we define 
  the {\it transition number} $n_{ik}(d)$ as
  the number of elements $b$ in $\{1,2,\ldots,(p-1)\}$ for which 
  $b\in
  G_i$ and $b+d \in G_k$. (Subscripts are taken modulo $3$. We do not
  count $b=p-d$).  Suppose now that $d \in G_a,$ i.e. that $d\equiv
  g^{a+3m}$ for some $m$ (congruences are modulo $p$).  For each $b$
  which contributes to $n_{ik}(1)$, we have $b\equiv g^{i+3u}$ and
  $b+1\equiv g^{k+3v}$ for some $u$ and $v$. Thus, from $d(b+1)=db+d$
  we get
\begin{equation}g^{k+a+3(m+v)}\,\equiv g^{i+a+3(m+u)}\,+\,d.
\end{equation}
Writing $n_{ik}$ instead of $n_{ik}(1)$, we thus get
\begin{equation} n_{i+a,k+a}(d)=n_{ik}. \label{eq:transnos}
\end{equation}
Let us now consider a simple cyclic $p$-root of index 3, and let
the corresponding $x$ be normalized by $x_{0}=1$ and have values given 
by (\ref{ind3}). Fix
$d$ such that $d \in G_a$, and consider the individual products in the
degree $d$ equation of (\ref{eq:cyk}). These products will take the
values $(c_k+a)/(c_i+a)$ with the frequency $n_{i+a,k+a}(d)$, the value
$c_a/1$ once (since $(p-1)\in G_0$), and the value $1/c_a$ once (since
$p-d \in G_a$). Thus (\ref{eq:transnos}) implies that all equations
whose
degrees $d$ belong to the same coset $G_a$, are identical, and the
system (\ref{eq:cyk})  
consists of the following 3 equations (where $n_{ik} = n_{ik}(1)$ are
the transition numbers, and the $c$ subscripts are counted modulo 3):
\begin{equation}
  {c_a \over 1}\hskip4pt+\hskip4pt{1\over
    c_{a}}\hskip4pt+\hskip4pt\sum_{k=0}^{2}\sum_{i=0}^{2}n_{ik}{c_{k+a}\over
    c_{i+a}}\hskip4pt=\hskip4pt 0,\hskip10pt (a=0,1,2).
\label{eqind3}\end{equation}

We will now return to the choice of the subscripts in $G_1$ and $G_2$.
Without loss of generality, we can (and do in fact from now on)
suppose that
\begin{equation}  n_{02}>n_{01}. 
\label{norder}\end{equation}
In fact, we must have $n_{02}\ne n_{01}$ (see Corollary \ref{ndiff}), and if
$n_{02}<n_{01}$, we replace the generator $g$ by $g':=g^{2+3j}$, for
some $j$ such that $2+3j$ is relatively prime to $p-1$. Since $g\in
G_1$ and $g'\in G_2$, this will interchange $G_1$ and $G_2$, and we
have arrived at (\ref{norder}).

  Finally, we will give the promised example: Let $p=13$, and take
  $g=2$ or 11. Then
  $G_0=\{1,5,8,12\},~G_1=\{2,3,10,11\},~G_2=\{4,6,7,9\},$ and we
  will have $n_{00}=0,~n_{01}=n_{10}=n_{12}=n_{21}=n_{22}=1$,~and
  $n_{02}=n_{20}=n_{11}=2$.
  
  \section{Number theoretic results used}
\label{sec-number}
\setcounter{equation}{0}

In this section we give some relations between the 
transition numbers $n_{ik}$ defined in (\ref {eq:transnos}) and 
appearing in (\ref {eqind3}). These relations will lead to explicit 
formulas for the $n_{ik}$.

The mapping $b\rightarrow p-b$
from $Z_{p}$ to $Z_{p}$ will leave each one of 
the sets $G_{i}$ invariant and thus we have 
\begin{equation}
n_{ij}=n_{ji},\,\,i,j=0,1,2.
\label{nsymm}
\end{equation} 
Moreover $\displaystyle \sum_{j=0}^2 n_{ij}=\sharp( 
(G_{i}\setminus \{p-1\})$, and thus (recall that we have defined 
$s={p-1\over 3}$)
\begin{equation}
\sum_{j=0}^2 n_{0j}= s-1,\,\,\,\, \sum_{j=0}^2 n_{1j}= 
\sum_{j=0}^2 n_{2j}= s.
\label{nrows}
\end{equation}
We will get one more linear relation between the $n_{ik}$ in the 
following way: By (\ref {eq:transnos}), all $n_{01}(d)$ with $d$ 
belonging to the same $G_{a}$ are equal. Thus, since 
$\sharp(G_{0})=\sharp(G_{1})=\sharp(G_{2})=s$, we get $\displaystyle
s\cdot s=\sum_{d=1}^{p-1}n_{01}(d)=\sum_{a=0}^2{s\cdot 
n_{-a,1-a}}$, which becomes 
\begin{equation} n_{01}+n_{12}+n_{20}=s.\label{ndiagonal}
\end{equation}
With the help of ({\ref{nsymm}), (\ref{nrows}) and (\ref{ndiagonal}) 
we can express
all our nine transition numbers $n_{ik}$ in terms of $n_{01}$ 
and $n_{02}$:

\begin{equation}
\label{nlist}
\left\{\begin{array}{lcl} n_{00} & =  & s-1-n_{01}-n_{02},\\
n_{11} & =  & n_{20} = n_{02},\\
 n_{22} & =  & n_{10} = n_{01},\\
  n_{12} & =  & n_{21} = s-n_{01}-n_{02}.\end{array}\right.
\end{equation}

These relations are given in \cite{pisa} and also in 
\cite{primesform}, Exercise 4.29 (d).
There is, however, one further equation satisfied by the transition 
numbers. We first state this equation in terms of $n_{12},n_{01}$ and
$n_{02}$:
\newtheorem  {n2eq}{Proposition}[section]
\begin{n2eq} Let $p$ be a prime $\equiv 1 \pmod{6}$, and let  
$n_{12},n_{01}$ and $n_{02}$ be the transition numbers defined in 
Section $\ref{sec-notation}$. Then
$$n_{01}n_{02}+n_{01}n_{12}+n_{02}n_{12}=n_{01}^2+n_{02}^2+n_{12}^2-n_{12}.$$
\label{n2eq}
\end{n2eq}
\vskip -8mm

We have proved this result by establishing 
the following explicit formulas for the 
convolutions $F*G$ (defined by $(F*G)(a)=\sum_{b\in{\bf Z_p}}F(a-b)G(b)$)
of certain complex-valued functions $F$ and $G$ on \Zp. Let $\Gamma_j$
be the characteristic functions $\chi_{G_j}$ of $G_j$ ($j=0,1,2)$,
and let $I=\chi_{\{0\}}$. Then, (with indices taken modulo 3):
$$\Gamma_i*\Gamma_i=n_{i,i}\Gamma_0+n_{i+2,i+2}\Gamma_1+n_{i+1,i+1}\Gamma_2+sI,$$
$$\Gamma_i*\Gamma_{i+1}=n_{i,i+1}\Gamma_0+n_{i+2,i}\Gamma_1+n_{i+1,i+2}\Gamma_2.$$
Our original proof of Proposition \ref{n2eq} used these formulas and the commutativity
and associativity of the convolution. Also, the reader of \cite{primesform}
is encouraged in Exercise 4.29 (e) to prove this proposition.
But it turns out that Proposition \ref{n2eq} is just a reformulation 
of a theorem of Gauss (in {\it Disquisitiones}, Article 358), 
which we give in a form a little more
precise than in \cite{IR} or \cite{silverman} or \cite{primesform}:
\newtheorem  {gauss}[n2eq]{Proposition}
\begin{gauss} Let $p$ be a prime $\equiv 1 \pmod{6}$, and let  
$n_{12},n_{01}$ and $n_{02}$ be the transition numbers defined in 
Section $\ref{sec-notation}$. Then there are integers $A$ and $B$ such that
$$4p=A^2+27B^2.$$
If we require that $A\equiv 1 \pmod{3}$ and $B>0$ (which is always 
possible and which we always do), then $A$ and $B$ are unique, and we have 
$$A=9n_{12}-p-1\,\, {\rm and}\,\, B =|n_{02}-n_{01}|.$$
\label{gauss}
\end{gauss}
\vskip -8mm
Since $4p$ is not a square, we must have $B\ne0$, and hence we get
the following corollary, which we needed at the end of Section 
\ref{sec-notation}:
\newtheorem {ndiff}[n2eq]{Corollary}
\begin{ndiff} Let $p$ be a prime $\equiv 1 \pmod{6}$, and let  
$n_{01}$ and $n_{02}$ be the transition numbers defined in 
Section $\ref{sec-notation}$. Then $n_{01}\ne n_{02}$.
\label{ndiff}
\end{ndiff}
Recall that we have in fact chosen $G_1$ and $G_2$ in such a way that 
$n_{02}>n_{01}$. Since $B>0$, we thus have 
\begin{equation}
\label{ABformulas}
A=9n_{12}-p-1\,\, {\rm and}\,\, B =n_{02}-n_{01}.
\end{equation}
 Solving the linear system given by (\ref{nlist}) and (\ref{ABformulas}) for 
 $n_{ik}$, we have proved the following corollary of Proposition 
 \ref{gauss}:
\newtheorem {nikformulas}[n2eq]{Corollary}
\begin{nikformulas} Let $p$ be a prime $\equiv 1 \pmod{6}$, let  
$n_{ik}$ be the transition numbers defined in 
Section $\ref{sec-notation}$, and let $A$ and $B$ be the numbers given in 
Proposition $\ref{gauss}$. Then 

\renewcommand{\arraystretch}{1.5}
\begin{equation}
\begin{array}{rcl}
n_{12} =n_{21}&=& \frac{1}{9}(p+A+1),\\
n_{02} =n_{20}=n_{11}&=& \frac{1}{18}(2p-A+9B-4),\\
n_{01}=n_{10}=n_{22}&=&\frac{1}{18}(2p-A-9B-4),\\
n_{00}+n_{11}+n_{22}&=&\frac{1}{3}(p-4).
\end{array}
\end{equation}

\label{nikformulas}
\end{nikformulas}

{\bf Proof of Proposition \ref{n2eq}:} Starting from Proposition \ref{gauss}
and replacing $A$ and $B$ by the expressions given there and then 
replacing $p$ by the expression  
$p=3(n_{01}+n_{12}+n_{20})+1$ from (\ref{ndiagonal}) we get
$$0=A^2+27B^2-4p=-36(n_{01}n_{02}+n_{01}n_{12}+n_{02}n_{12}
-n_{01}^2-n_{01}^2-n_{02}^2-n_{12}^2+n_{12})$$
which completes the proof.

{\bf Proof of Proposition \ref{gauss}:} The calculations needed are
given very explicitly in \cite{silverman}. In fact the theorem of Gauss
stated there in Section IV.2 is our Proposition \ref{gauss} except 
that
the {\it statement} of the theorem does not contain the value of $B$ and
for $A$ gives the value $M_p-p-1$, where $M_p$ is the number of 
solutions $(x,y,z)$ in ${\bf Z}_p^3$ of $x^3+y^3+z^3=0$ in the 
projective sense. In the {\it proof} of the theorem, the formula
$mB=[STT]-[STS]$ is given where $m$ is our $s$, where $R$ is our 
$G_0$, $S$ and $T$ are our $G_1$ and $G_2$ (in some order),
and where finally the symbol $[XYZ]$ is defined for subsets
$X,Y,Z$ of \Zp$\,$  as the number of triples $(x,y,z)$ such that
$x\in X$, $y\in Y$, and $z\in Z$ and $x+y+z=0$. In the course of the 
proof it is also shown that $mM_p=9[RTS]$. Thus all that remains for
us to have a proof of Proposition \ref{gauss} is to check that
$[G_1G_2G_2]-[G_1G_2G_1]=s(n_2-n_1)$ and $[G_0G_2G_1]=sn$.
We write $x+y+z=0$ as $x+y=-z$, and since $G_2=-G_2$, we have that
$$[G_{i+2}G_2G_{k+2}]=\sum_{y\in G_2}n_{i+2,k+2}(y)=sn_{ik},$$
where we have used (\ref{eq:transnos}) with $a=2$ and $d=y$.
Thus $[G_1G_2G_2]-[G_1G_2G_1]=s(n_{20}-n_{22})$ and 
$[G_0G_2G_1]=sn_{12}$,
and the result follows from (\ref{nlist}), which completes the proof.

\section{Reduction of the main problem}
\label{sec2}
\setcounter{equation}{0}
Let $p$ be a prime of the form $p=3s+1$, $s\in\bN$ and let
\[
4p = A^2+27B^2
\]
be the Gauss decomposition of $4p$, i.e.\ $A,B\in\bZ$, $A\equiv1 \pmod{3}$
and $B>0$ (cf. Proposition \ref{gauss}). Our main problem is to find 
all {\it simple} cyclic $p$-roots of index 3, i.e. to solve the set of 
equations
(cf. \ref{eqind3} and Corollary \ref{nikformulas})
\renewcommand{\arraystretch}{1.5}
\begin{equation}
\label{eq2-1}
\hskip8mm\left\{\begin{array}{lcl}
c_0+\frac{1}{c_0} &=&
-\frac{p-4}{3}-n_{12}\left(\frac{c_2}{c_1}+\frac{c_1}{c_2}\right)-n_{02}\left(\frac{c_0}{c_2}+
\frac{c_2}{c_0}\right)-n_{01}\left(\frac{c_1}{c_0}+\frac{c_0}{c_1}\right)\\
c_1+\frac{1}{c_1} &=&
-\frac{p-4}{3}-n_{12}\left(\frac{c_0}{c_2}+\frac{c_2}{c_0}\right)-n_{02}\left(\frac{c_1}{c_0}+
\frac{c_0}{c_1}\right)-n_{01}\left(\frac{c_2}{c_1}+\frac{c_1}{c_2}\right)\hskip14mm
\\
c_2+\frac{1}{c_2} &=&
-\frac{p-4}{3}-n_{12}\left(\frac{c_1}{c_0}+\frac{c_0}{c_1}\right)-n_{02}\left(\frac{c_2}{c_1}+
\frac{c_1}{c_2}\right)-n_{01}\left(\frac{c_0}{c_2}+\frac{c_2}{c_0}\right)\end{array}\right.
\end{equation}
with
\begin{equation}
\label{eq2-2}
n_{12} = \frac{p+A+1}{9},\quad n_{02} = \frac{2p-A+9B-4}{18},\quad 
n_{01}=\frac{2p-A-9B-4}{18}.
\end{equation}

\begin{proposition}
\label{prop2-1}
Assume $(c_0,c_1,c_2)$ is a solution to  {\rm\eqref{eq2-1}}. Then the 
numbers
\begin{equation}
\label{hdef}
h_j = \frac{c_{j+2}}{c_{j+1}}+\frac{c_{j+1}}{c_{j+2}},\quad j=0,1,2,
\end{equation}
(index counted modulo {\rm 3}) are up to a cyclic permutation given by
\begin{equation}
\label{eq2-3}
h_j = \xi_1+\eta_1 \cos\left(\theta -\frac{2\pi}{3}j\right),\quad j= 
0,1,2,
\end{equation}
where $\theta=\frac13 \Arccos\left(\frac{A}{2\sqrt{p}}\right)$ and the
pair $(\xi_1,\eta_1)$ is one of the following 4 pairs:

\renewcommand{\arraystretch}{1.5}
\begin{equation}
\label{eq2-4}
\left\{ \begin{array}{lcl} \xi_1^{(0)} &=& 2\\
\eta_1^{(0)} &=& 0\,, \end{array}\right.
\end{equation}
{\large\begin{equation}
\label{eq2-5}
\left\{ \begin{array}{lcl} \xi_1^{(1)} &=& 
-\frac{p^2-6p+2A}{p^2-3p-A}\\
\eta_1^{(1)} &=& \frac{6\sqrt{p}(p-4)}{p^2-3p-A}\,, \end{array}\right.
\end{equation}

\begin{equation}
\label{eq2-6}
\left\{\begin{array}{lcl} \xi_1^{(2)} &=& 
\frac{-2pA-9p-4+3\sqrt{p(p+4A+16)}}{2(pA+3p-1)}\\
\eta_1^{(2)} &=& \frac{3\sqrt{p}(p+2)-3p\sqrt{p+4A+16}}{pA+3p-1}\,,
\end{array}\right.
\end{equation}

\begin{equation}
\label{eq2-7}
\left\{\begin{array}{lcl} \xi_1^{(3)} &=& 
\frac{-2pA-9p-4-3\sqrt{p(p+4A+16)}}{2(pA+3p-1)}\\
\eta_1^{(3)} &=& \frac{3\sqrt{p}(p+2)+3p\sqrt{p+4A+16}}{pA+3p-1}\,.
\end{array}\right.
\end{equation}}
\renewcommand{\arraystretch}{1.0}
\end{proposition}

\begin{remark}
\label{rem2-2}
{\rm 
a) Let us first check that all the above formulas give
well-defined real numbers:
Since $p>4$ and $|A|<2\sqrt{p}$ we have
\[
p^2-3p-A > p^2-3p-2\sqrt{p} = \sqrt{p}(\sqrt{p}-2)(\sqrt{p}+1)^2 >0.
\]
Moreover,
\[
p+4A+16 > p-8\sqrt{p}+16 = (\sqrt{p}-4)^2\ge 0
\]
and since $A\equiv1 \pmod{3}$, we have $|A+3|\ge 1$. Hence
\[
|pA+3p-1|\ge |(A+3)p|-1\ge p-1 > 0.
\]

b)We do not prove in this section that all four cases
\eqref{eq2-4}--\eqref{eq2-7} actually occur. However this will follow
from the proof of Theorem \ref{thm3-1} in the next section.
}
\end{remark}

\noindent{\bf Proof of Proposition \ref{prop2-1}:} \  To make our 
method
of proof more transparent, we first consider the case $p=7$. In this
case $A=B=1$, $n_{12}=n_{02}=1$, and $n_{01}=0$. Put
\[
f_j = c_j+\frac{1}{c_j}\qquad\mbox{and}\qquad h_j =
\frac{c_{j+2}}{c_{j+1}}+\frac{c_{j+1}}{c_{j+2}}.
\]
Then \eqref{eq2-1} becomes
\begin{equation}
\label{eq2-9}
\left\{\begin{array}{lcl} f_0 &=& -1-h_0-h_1\\
f_1 &=& -1-h_1-h_2\\
f_2 &=& -1-h_2-h_0 .\end{array}\right.
\end{equation}
Consider now the matrix
\[
K = \left[\begin{array}{cccc}
2&f_0&f_1&f_2\\
f_0&2&h_2&h_1\\
f_1&h_2&2&h_0\\
f_2&h_1&h_0&2\end{array}\right].
\]  
Since
\[
K =\left[\begin{array}{cccc}1\\c_{0}\\c_{1}\\c_{2}\end{array}\right]
[1,\frac{1}{c_0},\frac{1}{c_1},\frac{1}{c_2}]+
\left[\begin{array}{cccc}1\\\frac{1}{c_{0}}\\\frac{1}{c_{1}}\\
\frac{1}{c_{2}}\end{array}\right][1,c_0,c_1,c_2],\]
we get (considering $K$ as an operator on column vectors)
\[
\range(K)=\Span\left\{\left[\begin{array}{cccc}1\\c_{0}\\c_{1}\\c_{2}\end{array}
\right],
\left[\begin{array}{cccc}1\\\frac{1}{c_{0}}\\\frac{1}{c_{1}}\\
\frac{1}{c_{2}}\end{array}\right]\right\}.\]
Hence $\rank(K)\le 2$, and thus all $3\times 3$ submatrices of $K$ 
have
determinant $=0$.

Let $L=(\ell_{ij})^4_{i,j=1}$ be the co-factor matrix of $K$, i.e.
\[
\ell_{ij}=(-1)^{i+j}\det(K_{ij}),
\]
where $K_{ij}$ is the $3\times 3$ minor of $K$ obtained by erasing the
$i$'th row and the $j$'th column. Put
\begin{equation}
\label{eq2-p}
\left\{\begin{array}{lcl}
p_1 &=& \ell_{11}\\
p_2 &=& \ell_{12}+\ell_{13}+\ell_{14}\\
p_3 &=& \ell_{22}+\ell_{33}+\ell_{44}\\
p_4 &=& \ell_{23}+\ell_{34}+\ell_{42}.\end{array}\right.
\end{equation} 
Since $\ell_{ij}=0$ for all $i$ and $j$, we have in particular
\[
p_1=p_2=p_3=p_4=0.
\]
This gives four equations of degree three in $(f_0,f_1,f_2,h_0,h_1,h_2)$, but
taking \eqref{eq2-9} into account, we can consider $p_1,p_2,p_3,p_4$ 
as
polynomials in $(h_0,h_1,h_2)$ only, namely
\begin{eqnarray*}
p_1 &=& 8-2(h_0^2+h_1^2+h_2^2)+2h_0h_1h_2\\
p_2 &=& 12-4(h_0+h_1+h_2)-3(h_0^2+h_1^2+h^2_2)
-4(h_0h_1+h_1h_2+h_2h_0)\\
&& \quad 
-(h^3_0+h^3_1+h^3_2)+2(h_0h_1^2+h_1h^2_2+h_2h_0^2)+3h_0h_1h_2,\\
p_3 &=& 
12-14(h_0+h_1+h_2)-8(h^2_0h^2_1+h^2_2)-2(h_0h_1+h_1h_2+h_2h_0)\\
&& \quad +2(h_0h_1^2+h_1h_2^2+h_2h_0^2)+4(h_0^2h_1+h_1^2h_2+h_2^2h_0) 
6h_0h_1h_2,\\
p_4 &=& 6+3(h_0+h_1+h_2)+(h_0^2+h_1^2+h_2^2)+5(h_0h_1+h_1h_2+h_2h_0)\\
&& \quad -2(h_0h_1^2+h_1h_2^2+h_2h_0^2)-6h_0h_1h_2.
\end{eqnarray*}
Let $s_1,s_2,s_3$ denote the three elementary symmetric polynomials in
$h_0,h_1,h_2$:
\begin{equation}
\label{eq2-10}
\left\{\begin{array}{lcl} s_1 &=& h_0+h_1+h_2\\
s_2&=& h_0h_1+h_1h_2+h_2h_0\\
s_3&=& h_0h_1h_2\end{array}\right.
\end{equation}
and let $a$ denote the antisymmetric polynomial:
\begin{equation}
\label{eq2-11}
a=(h_0-h_1)(h_1-h_2)(h_2-h_0).
\end{equation}
Then,
\begin{eqnarray*}
h_0^2+h_1^2+h_2^2 &=& s_1^2-2s_2\\
h_0^3+h_1^3+h_2^3&=& s_1^3-3s_1s_2+3s_3\\
h_0h_1^2+h_1h_2^2+h_2h_0^2&=& \tfrac12(s_1s_2-3s_3+a)\\
h_0^2h_1+h_1^2h_2+h_2^2h_0 &=& \tfrac12(s_1s_2-3s_3-a).
\end{eqnarray*}
Hence $p_1,p_2,p_3,p_4$ can be expressed as polynomials in 
$s_1,s_2,s_3$
and $a$. One gets
\begin{eqnarray*}
p_1 &=& (8-2s_1^2)+4s_2+2s_3\\
p_2 &=& (12+4s_1-3s_1^2-s_1^3)+(2+4s_1)s_2-3s_3+a\\
p_3 &=& (12-14s_1-8s_1^2)+(14+3s_1)s_2-3s_3-a\\
p_4 &=& (6+3s_1+s_1^2)+(3-s_1)s_2-3s_3-a.
\end{eqnarray*}
Therefore the equations $p_1=p_2=p_3=p_4=0$ can be rewritten in the 
form
\begin{equation}
\label{eq2-13}
\left[\begin{array}{ccrr} 8-2s_1^2 & 4 & 2 & 0\\ 12+4s_1-3s_1^2-s_1^3 
&
2+4s_1 & -3 & 1\\
12-14s_1-8s_1^2 & 14+3s_1 & -3 & -1\\ 6+3s_1+s_1^2 & 3-s_1 & -3 &
-1\end{array}\right] \left[\begin{array}{c}
1\\s_2\\s_3\\a\end{array}\right] = \left[\begin{array}{c}
0\\0\\0\\0\end{array}\right] .
\end{equation}
A necessary condition for the existence of solutions to this system of
equations is that the determinant of the coefficient matrix $M$ is
0. One finds
\[
\det(M) = 8(s_1-6)(s_1+1)(s_1^2+9s_1+15).
\]
Thus $s_1$ must be one of the 4 numbers
\[
s_1^{(0)}=6,\quad s_1^{(1)}=-1,\quad s_1^{(2)}=\frac{-9+\sqrt{21}}{2}\quad\mbox{or}\quad
s_1^{(3)}=\frac{-9-\sqrt{21}}{2}.
\]
Let $M^{(i)}$ be the matrix obtained by substituting $s_1=s_1^{(i)}$ 
in
$M$ ($i=0,1,2,3$). It is easy to compute the kernel for $M^{(i)}$,
$i=0,1,2,3$. One finds $\dim(\ker(M^{(i)})=1$ in all cases, and (for 
convenience writing vectors in row form)
\renewcommand{\arraystretch}{1.5}
\begin{eqnarray*}
\ker(M^{(0)}) &=& \Span\big\{[1,12,8,0]\big\}\\
\ker(M^{(1)}) &=& \Span\big\{[1,-2,1,7]\big\}\\
\ker(M^{(2)}) &=&
\Span\bigg\{\bigg[1,-9+\sqrt{21},\frac{79-17\sqrt{21}}{2},-189+42\sqrt{21}\bigg]\bigg\}\\
\ker(M^{(3)}) &=& 
\Span\bigg\{\bigg[1,-9-\sqrt{21},\frac{79+17\sqrt{21}}{2},
-189+42\sqrt{21}\bigg]\bigg\}. 
\end{eqnarray*}
Hence there are exactly 4 solutions ($s_1,s_2,s_3,a$) to 
\eqref{eq2-13}:

\begin{equation}
\label{sa7}
\begin{array}{lllll}
(s_1^{(0)}= 6,& s_2^{(0)}=12,&s_3^{(0)}=8,&a^{(0)}=0&),\\
(s_1^{(1)}=-1,&s_2^{(1)}=-2,&s_3^{(1)}=1,&a^{(1)}=-7&),\\
\bigg(s_1^{(2)}=\frac{-9+\sqrt{21}}{2},& 
s_2^{(2)}=-9+2\sqrt{21},&
s_3^{(2)}=\frac{79-17\sqrt{21}}{2},& 
a^{(2)}=-189+42\sqrt{21}&\bigg),\\
\bigg(s_1^{(3)}=\frac{-9-\sqrt{21}}{2},& 
s_2^{(2)}=-9-2\sqrt{21},&
s_3^{(2)}=\frac{79+17\sqrt{21}}{2},& 
a^{(3)}=-189-42\sqrt{21}&\bigg).
\end{array}
\end{equation}

However, there is a hidden relation between $s_1,s_2,s_3$ and $a$, namely
$a^2$ is a symmetric polynomial in $(h_0,h_1,h_2)$ and can therefore 
be
expressed in terms of $s_1,s_2$ and $s_2$. One finds
\begin{equation}
\label{eq2-14}
a^2=s_1^2s_2^2-4s_1^3s_3-4s_2^3+18s_1s_2s_3-27s_3^2.
\end{equation}
It is elementary to check that this equality holds for each of the 
four
sets $(s_1^{(i)},s_2^{(i)},s_3^{(i)},a^{(i)})$ found above.

We must now in each case find $h_0,h_1,h_2$ by solving the 4 
equations:
\begin{equation}
\label{eq2-15}
\left\{\begin{array}{l} h_0+h_1+h_2 = s_1^{(i)}\\
h_{0}h_{1}+h_{1}h_{2}+h_{2}h_{0}=s_{2}^{(i)}\\
h_0h_1h_2 = s_3^{(i)}\\
(h_0-h_1)(h_1-h_2)(h_2-h_0) = a^{(i)}.\end{array}\right .
\end{equation}
The solutions $(h_0,h_1,h_2)$ to the first 3 equations in 
\eqref{eq2-15}
are exactly the three roots (in arbitrary order) to the polynomial
\begin{equation}
\label{eq2-16}
h^3-s_1^{(i)}h^2+s_2^{(i)}h-s_3^{(i)}.
\end{equation}
Since \eqref{eq2-14} holds in each of the four cases $i=0,1,2,3$, we 
have
\[
(h_0-h_1)(h_1-h_2)(h_2-h_0) = \pm a^{(i)}.
\]
Hence the 4'th coordinate in the solution to the equations \eqref{eq2-13} only determines the cyclic
order of the three numbers $(h_0,h_1,h_2)$. For $i=0$, \eqref{eq2-16}
becomes
\[
h^3-6h^2+12h-8=0.
\]
Hence $h_0=h_1=h_2=2$ which corresponds to case \eqref{eq2-4} in
Proposition \ref{prop2-1}.

In the cases $i=1,2,3$ we solve \eqref{eq2-16} by the classical 
trigonometric formula in
the form of Lemma \ref{lemma2-4} below, where we use
\eqref{eq2-25} when $a<0$ and \eqref{eq2-26} when $a>0$. This will give
the correct cyclic order of $(h_0,h_1,h_2)$. Note that Lemma
\ref{lemma2-4} can be applied because in all 3 cases $(i=1,2,3)$
$s_1,s_2,s_3$ and $a$ are all real (being solutions to the real 
linear system \eqref{eq2-13}) and thus $a^2>0$, which by 
\eqref{eq2-14} means that 
$s_1^2s_2^2-4s_1^3s_3-4s_2^3+18s_1s_2s_3-27s_2^3=a^2>0$. Hence, up to
cyclic permutation of $(h_0,h_1,h_2)$ we have
\[
h_j= \xi_1+\eta_1\cos\bigg(\theta -\frac{2\pi}{3}j\bigg)
\]
where
\renewcommand{\arraystretch}{1.5}
{\large\[
\left\{\begin{array}{lcl} \xi_1 &=& \frac13 s_1\\
\eta_1 &=& -\sign(a)\cdot\frac23(s_1^2-3s_2)^{\frac12}\\ \theta&=&
\frac13\Arccos\big(-\sign(a)\frac{2s_1^3-9s_1s_2+27s_3 
}{2(s_1^2-3s_2)^{\frac32}}\big) \end{array}
\right. 
\].}
\renewcommand{\arraystretch}{1.0}
It turns out that
${\displaystyle\theta^{(i)}=\frac13\Arccos\left(\frac{1}{2\sqrt{7}}\right)}$ in all
three cases $(i=2,3,4)$, while
\begin{eqnarray*}
\left(\xi_1^{(1)},\eta_1^{(1)}\right) &=&
\left(-\frac13,\frac23\sqrt{7}\right)\\
\left(\xi_1^{(2)},\eta_1^{(2)}\right) &=&
\left(-\frac32+\frac{\sqrt{21}}{6},\sqrt{7}-\frac73\sqrt{3}\right)\\
\left(\xi_1^{(3)},\eta_1^{(3)}\right) &=&
\left(-\frac32-\frac{\sqrt{21}}{6},\sqrt{7}+\frac73\sqrt{3}\right).
\end{eqnarray*}
This gives case \eqref{eq2-5}, \eqref{eq2-6}, and \eqref{eq2-7} 
respectively
in Proposition \ref{prop2-1} in the case $p=7$.

Consider now a general prime $p$, $p\equiv1 \pmod{3}$. This case is
mathematically no more difficult than the case $p=7$ but a
computer algebra language as MAPLE or MATHEMATICA is helpful for bookkeeping
purpose. Using \eqref{eq2-2} and \eqref{eq2-1} instead of \eqref{eq2-9}, the polynomials
\eqref{eq2-10} again becomes polynomials in $s_1,s_2,s_3,a$, namely
\begin{equation}
\label{eq2-17}
 \left[\begin{array}{c} p_1\\p_2\\p_3\\p_4\end{array}\right] =
 \left[\begin{array}{cccc} m_{11} & m_{12} & m_{13} & m_{14}\\ m_{21} 
&
 m_{22} & m_{23} & m_{24}\\ m_{31} & m_{32} & m_{33} & m_{34}\\ 
m_{41} &
 m_{42} & m_{43} & m_{44}\end{array}\right] \left[\begin{array}{c} 
1\\s_2\\s_3\\a\end{array}\right]
\end{equation}
where the $m_{ij}$:s are the following 16 polynomials in $s_1$:
\renewcommand{\arraystretch}{1.5}
\begin{eqnarray*}
m_{11} &=& -2s_1^2+8\\
m_{12} &=& 4\\
m_{13} &=& 2\\
m_{14} &=& 0\\
m_{21} &=& \tfrac19(A+p+1)s_1^3+\tfrac19(2A-7p+20)s_1^2+4s_1+(4p-16)\\
m_{22} &=& \tfrac49(A+p+1)s_1+\tfrac13(4p-2A-20)\\
m_{23} &=& -A-2\\
m_{24} &=& B\\
m_{31} &=&
\tfrac{2}{81}(p^2-pA-7p+a^2+2A+A)s_1^3-\tfrac{2}{27}(pA+12p+17)s_1^2\\
&&\quad -\tfrac23p(p-4)s_1+\tfrac43(2-p^2+8p)\\
m_{32} &=& \tfrac{1}{27}(-6A-12p-3A^3-8+2pA)s_1 + \tfrac29(6p+pA+14)\\
m_{33} &=& \tfrac13(2A+A^2-2p+2)\\
m_{34} &=& -\tfrac13(A+2)B\\
m_{41} &=&
\tfrac18(7A-p^2+4p+2A^2+pA+5)s_1^3+\tfrac{1}{27}(-6A+6p+pA-16) s_1\\
&&\quad+\tfrac13(p^2-4p-12)s_1+\tfrac23(p^2-8p+16)\\
m_{42} &=& -\tfrac{1}{27}(9A+3A^2+pA+8)s_1+\tfrac19(6A-pA+28)\\
m_{43} &=& \tfrac13(2A+A^2-2p+2)\\
m_{44} &=& -\tfrac13(A+2)B.
\end{eqnarray*}
\renewcommand{\arraystretch}{1.0}
Since $p_1=p_2=p_3=p_4=0$, we must have $\det M=0$ where
$M=(m_{ij})_{i,j=1}^4$. One finds
\[
\det M=\tfrac{8B}{729}(s_1-6)q(s_1)r(s_1),
\]
where
\begin{equation}
\label{eq2-18}
\left\{\begin{array}{lcl} q(s_1) &=& (p^2-3p-A)s_1+(6A+3p^2-18p)\\
r(s_1) &=&
(pA+3p-1)s_1^2+(6pA+27p+12)s_1+(9pA+54p-36).\end{array}\right. 
\end{equation}
It is interesting that if $\det M$ is considered as a polynomial in 
the independent variables $s_{1},p,A,B$, forgetting the relation $4p = 
A^2+27B^2$, we will get an irreducible cubic polynomial instead of 
$q(s_{1})r(s_{1})$.
By Remark \ref{rem2-2}, $p^2-3p-A\ne 0$ and $pA+3p-1\ne 0$, so the
equation $\det(M)=0$ has exactly 4 solutions (counted with
multiplicity), namely
\renewcommand{\arraystretch}{1.5}
\begin{equation}
\label{eq2-19}
\left\{\begin{array}{lcl} s_1^{(0)} &=& 6\\ s_1^{(1)} &=&
\frac{18p-3p^2-6A}{p^2-3p-A}\\
s_1^{(2)} &=& \frac{-6pA-27p-12+9\sqrt{p(p+4A+16)}}{2(pA+3p-1)}\\
s_1^{(3)} &=& \frac{-6pA-27p-12-9\sqrt{p(p+4A+16)}}{2(pA+3p-1)}.
\end{array}\right. 
\end{equation}
\renewcommand{\arraystretch}{1.0}
Let $M^{(i)}$ be the $4\times 4$-matrix  obtained by substituting
$s_1=s_1^{(i)}$ in $M$. We next compute the kernel for $M^{(i)}$ in each 
of the four cases. Let $M_{jk}^{(i)}$ be the $3\times 3$ minor of 
$M^{(i)}$ obtained by erasing the
$j$'th row and the $k$'th column of $M^{(i)}$. Then
\[
\det(M_{11}^{(i)}) =
-\frac{2B}{27}(A+p+1)((pA+A+4p)s_1^{(i)}+3pA-6A+12p).
\]
In particular
\begin{eqnarray*}
\det(M_{11}^{(0)}) &=& -\frac{2B}{3} p(p+A+1)(A+4),\\
\det(M^{(1)}_{11}) &=& -\frac{2B}{3}p(p+A+1)(4p-A^2),\\
\det(M_{11}^{(2)})\cdot\det(M_{11}^{(3)}) &=&
-\frac{4B^2(p+A+1)^2(A+4)(4p-A^2)}{9(pA+3p-1)}.
\end{eqnarray*}
Since $A\equiv1\pmod 3$, we have $A+4\ne 0$. Moreover $4p-A^2=27B^2>0$ and
$p+A+1>(\sqrt{p}-1)^2\ge 0$. Hence $\det(M_{11}^{(i)})\ne 0$ in all 4
cases. Together with $\det(M^{(i)})=0$, this shows that for all $i,s$,
$M^{(i)}$ has rank 3 and thus
\[
\dim(\Ker(M^{(i)})=1,\qquad i=0,1,2,3.
\]
Hence in each case $(i=0,1,2,3)$, $s_2^{(i)},s_3^{(i)}$ and $a^{(i)}$
are uniquely determined by \eqref{eq2-17}. Applying Cramer's rule to 
the
last three equations in \eqref{eq2-17} we get
\[
s_2^{(i)} = -\frac{\det M_{12}^{(i)}}{\det M_{11}^{(i)}},\quad
s_3^{(i)}=\frac{\det M_{13}^{(i)}}{\det M_{11}^{(i)}},\quad
a^{(i)}=-\frac{\det M_{14}^{(i)}}{\det M_{11}^{(i)}}.
\]
For $i=0$, $(s_1^{(0)},s_2^{(0)},s_3^{(0)},a^{(0)})=(6,12,8,0)$ as in
the case $p=7$ and for $i=1$ we have
\renewcommand{\arraystretch}{1.5}
\begin{equation}
\label{sssa1}
\left\{\begin{array}{lcl} s_1^{(1)} &=& \frac{18p-3p^2-6A}{p^2-3p-A}\\
s_2^{(1)} &=&
3\frac{4p^2A-24pA+4A^2+p^4-21p^3+108p^2-144p}{(p^2-3p-A)^2}\\
s_3^{(1)} &=& 
\frac{20p^2A-96pA+8A^2-p^4+4p^3-360p^2+864p}{(p^2-3p-A)^2}
\\
a^{(1)} &=& -\frac{729p(p-4)^3 B}{p^2-3p-A)^3} .\end{array} \right.
\end{equation}
\renewcommand{\arraystretch}{1.0}

For  $i=2,3$, it is more convenient to express the solutions in terms 
of
$u=\sqrt{p}$ and $v=\sqrt{p+4A+16}$. We get
\renewcommand{\arraystretch}{1.5}

\begin{equation}
\label{sssa2}
\left\{\begin{array}{lcl} s_1^{(2)} &=& -3\frac{u^2+uv-4}{u^2+uv+2}\\
s_2^{(2)} &=&
3\frac{(u^2+uv+6u-4)(u^2+uv-6u-4)}{(u^2+uv+2)^2}\\
s_3^{(2)} &=& \frac{u^4+2u^3v-176u^2+u^2v^2+40uv-32}{(u^2+uv+2)^2}
\\
a^{(2)} &=& 5832 \frac{Bu^2}{(u^2+uv+2)^3}\end{array} \right.
\end{equation}
and
\begin{equation}
\label{sssa3}
\left\{\begin{array}{lcl} s_1^{(3)} &=& -3\frac{u^2-uv-4}{u^2-uv+2}\\
s_2^{(3)} &=&
3\frac{(u^2-uv+6u-4)(u^2-uv-6u-4)}{(u^2-uv+2)^2}\\
s_3^{(3)} &=& \frac{u^4-2u^3v-176u^2+v^2u^2-40uv-32}{(u^2-uv+2)^2}
\\
a^{(3)} &=& 5832 \frac{Bu^2}{(u^2-uv+2)^3} .\end{array} \right.
\end{equation}
\renewcommand{\arraystretch}{1.0}
Note that all the numbers are well-defined because by Remark
\ref{rem2-2}, $p^2-3p-A >0$, $p+4A+16>0$ and
\[
(u^2+uv+2)(u^2-uv+2) = -4(pA+3p-1)\ne 0.\]

For $i=0$, we get as for $p=7$ that $h_0=h_1=h_2=2$ which corresponds 
to
\eqref{eq2-4} in Proposition \ref{prop2-1}. It is easy to check that 
the
identity \eqref{eq2-14} is satisfied for the above sets
$(s_1^{(i)},s_2^{(i)},s_3^{(i)},a^{(i)})$, so as in the case $p=7$ we
can determine $h_0,h_1,h_2$ by Lemma \ref{lemma2-4} where we use
\eqref{eq2-25} when $a^{(i)}<0$ and \eqref{eq2-26}, when $a^{(i)}>0$ 
to
obtain the correct cyclic ordering. Note that $a^{(1)}<0$, $a^{(2)}>0$
and $\sign(a^{(3)})=\sign(u^2-uv+2)=-\sign(pA+3p-1)$. We obtain
\[
h_j = \xi^{(i)}_1+\eta^{(i)}_1
\cos\left(3\theta^{(i)}-\frac{2\pi}{3}j\right),\quad j=0,1,2,
\]
where $\theta^{(i)}=\frac13\Arccos\left(\frac{A}{2\sqrt{p}}\right)$ in
all three cases $(i=2,3,4)$, while
\renewcommand{\arraystretch}{1.5}
\begin{equation}
\label{eq2-20}
\left\{\begin{array}{lcl}
(\xi_1^{(1)},\eta_1^{(1)}) &=& \left(-\frac{p^2-6p+2A}{p^2-3p-A},
\frac{6\sqrt{p}(p-4)}{p^2-3p-A}\right)\\
(\xi_1^{(2)},\eta_1^{(2)}) &=& \left(-\frac{u^2+uv-4}{u^2+uv+2},
-\frac{12u}{u^2+uv+2}\right)\\
(\xi_1^{(3)},\eta_1^{(3)}) &=& \left(-\frac{u^2-uv-4}{u^2-uv+2},
-\frac{12u}{u^2-uv+2}\right).\end{array}\right.
\end{equation}
\renewcommand{\arraystretch}{1.0}
Using $u=\sqrt{p}$, $v=\sqrt{p+4A+16}$, we get \eqref{eq2-5},
\eqref{eq2-6} and \eqref{eq2-7} in Proposition \ref{prop2-1}.
This completes the proof of Proposition \ref{prop2-1},
\begin{remark}{\rm
\label{anot0}
It easily follows from the proof that if $c=(c_{0}, c_{1}, c_{2})$ is a solution  to
the system \eqref{eq2-1} and two $c_{i}$ are equal, then they are all 
equal. In fact, if e.g. $c_{1}= c_{2}$, then with $h$ as in 
\eqref{hdef} we get $h_{1}= h_{2}$, which leads to $a=0$. But since $B\ne 0$,
it follows from \eqref{sssa1}, \eqref{sssa2}, 
and \eqref{sssa3} that $a\ne 0$ in all cases except the case 
where all $h_{i}=2$.
}
\end{remark}
\begin{remark}{\rm
\label{rem2-4}
(a) At a first glance it is surprising that the angle $\theta$ 
in
the solution formula above is the same for $i=1,2,3$. However, this 
fact
has a fairly simple explanation: Computing the linear combination
\[
(p-1)p_1-2p_2-p_3-2p_4
\]
of the polynomials $p_i=p_i(s_1,s_2,s_3,a)$ given by \eqref{eq2-17} one gets
\[
\frac{4p-A^2}{27}(2s^3_1-9s_1s_2+27s_3)+ABa.
\]
Since $p_1=p_2=p_3=p_4=0$ and $B^2=\frac{4p-A^2}{27}$, we have the
following identity
\begin{equation}
\label{eq2-21}
B(2s_1^3-9s_1s_2+27s_3)+Aa=0.
\end{equation}
But if $h_j=\xi_1+\eta_1\cos\left(\theta-\frac{2\pi}{3} j\right)$,
$j=0,1,2$, and $s_1,s_2,s_3,a$ are defined as in \eqref{eq2-10} and
\eqref{eq2-11} one finds
\[
2s_1^3-9s_1s_2+27s_3=\frac{27}{4}\eta^3_1\cos 3\theta
\]
and
\[
a=-\frac{3\sqrt{3}}{4}\eta_1^3\sin 3\theta.
\]
Hence, when $\eta_{1}\ne 0$, \eqref{eq2-21} is equivalent to
\[
3\sqrt{3} B\cos 3\theta -A\sin 3\theta =0.
\]
This has a unique solution $\theta\in\left(0,\frac{\pi}{3}\right)$,
namely
\[
\theta = \frac13\Arccot\left(\frac{A}{3\sqrt{3} B}\right) = \frac13
\Arccos\left(\frac{A}{2\sqrt{p}}\right).
\]
(b) It is interesting to compare the solutions in Proposition
\ref{prop2-1} with the Gaussian cubic sum
\[
G = \sum^{p-1}_{j=0} e^i \frac{2\pi}{p} j^3.
\]
It is known that (cf.~\cite{heath} or Section IV.2 of \cite{silverman}) that for $p$ 
prime,
$p\equiv1 \pmod{3}$, $G$ is a solution to the cubic equation
\[
x^3-3px-pA=0.
\]
This equation has the 3 solutions
\[
x_j = 2\sqrt{p}\cos\left(\theta-\frac{2\pi}{3} j\right),\qquad j=0,1,2
\]
where $\theta=\frac13 \Arccos\left(\frac{A}{2\sqrt{p}}\right)$ as in
Proposition \ref{prop2-1}.
}
\end{remark}

It is a famous problem (the Problem of Kummer) to decide for each $p$
which of the 3 solutions is equal to $G$ (cf.\cite{heath} and Section 9.12 of \cite {IR} 
or Section IV.2 of \cite{silverman}). 
 
We conclude this section by stating as a lemma the classical 
trigonometric solution of a cubic equation with 3 real roots. For 
completeness, we recall an elementary proof (cf e.g. \S 47 of \cite{dickson})

\begin{lemma} \label{lemma2-4} 
Consider the cubic equation
\begin{equation}
\label{eq2-22}
h^3-s_1h^2+s_2h-s_3=0
\end{equation}
and assume that $s_1,s_2,s_3\in\bR$ and
\begin{equation}
\label{eqa2>0}
s_1^2s_2^2-4s_1^3s_3-4s^3_2+18s_1s_2s_3-27s_3^2 > 0.
\end{equation}
Then
\begin{eqnarray}
\label{eq2-23}
s_1^2-3s_2 &>& 0,\\
\label{eq2-24}
|2s_1^3-9s_1s_2+27s_3| &<& 2(s_1^2-3s_2)^\frac32.
\end{eqnarray}
Moreover \eqref{eq2-22} has {\rm 3} different real solutions. Listed in
decreasing order $h_0>h_1>h_2$, the solutions are
\begin{equation}
\label{eq2-25}
h_j = \frac{s_1}{3} + \frac23 (s_1^2-3s_2)^\frac12
\cos\left(\theta-\frac{2\pi j}{3}\right),\quad j=0,1,2,
\end{equation}
where 
\begin{equation}
\label{eq2-theta}
\theta = \frac13
\Arccos\left(\frac{2s_1^3-9s_1s_2+27s_3}{2(s_1^2-3s_2)^\frac32}\right).
\end{equation}
and listed in increasing order $h'_0<h'_1<h'_2$, the solutions are
\begin{equation}
\label{eq2-26}
h'_j = \frac{s_1}{3} + \frac23 (s_1^2-3s_2)^\frac12
\cos\left(\theta'-\frac{2\pi j}{3}\right),\quad j=0,1,2
\end{equation}

where 
\begin{equation}
\label{eq2-thetaprim}
\theta' = 
\tfrac13{\Arccos\left(-\frac{2s_1^3-9s_1s_2+27s_3}{2(s_1^2-3s_2)^\frac32}}\right).
\end{equation}

\end{lemma}

\bproof Let $h_{0},h_{1},h_{2}$ be the solutions to \eqref{eq2-22} 
and define  $a$ by \eqref{eq2-11}. Then from \eqref{eq2-14} follows that 
(the discriminant) $a^2>0$. Hence
$h_{0},h_{1},h_{2}$ are real and different (since if e.g. $h_{1}= c+i 
d, h_{2}=c- i d$ with $d\ne 0$ and $h_{0}\in\bR$ 
we would have
$a^2=-4d^2((h_{0}-c)^2+d^2)^2<0$, whereas e.g. $h_{0}=h_{1}$ would imply 
that $a=0$
).
Substituting $h=u+\frac{s_1}{3}$ in equation \eqref{eq2-22} we get 
\begin{equation}
\label{eq2-27}
u^3+ru+q=0
\end{equation}
where $r=-\frac13 (s_1^2-3s_2)$ and
$q=-\frac{1}{27}(2s_1^3-9s_1s_2+27s_3)$. Applying \eqref{eq2-14} with 
$s_{1}=0,s_{2}=r,s_{3}=-q$ we get $a^2=-4r^3-27q^2$.
Since the transformation from $h$ to $u$ is a translation, the 
discriminant does not change and thus \eqref{eqa2>0} becomes 
$-4r^3-27q^2>0$. Thus $r<0$, which is \eqref{eq2-23}.
Next we consider \eqref{eq2-24}. Squaring this relation and 
introducing $r$ and $q$ we give it the form $|-27q|^2< 4(-3r)^3$, 
which we have just seen is true.

Taking $u=m z$ in \eqref{eq2-27} we get the 
equation
\begin{equation}
\label{eq2-uz}
z^3+\frac{r}{m^2} z +\frac{q}{m^3}=0.
\end{equation}

We now start from the trigonometric identity 
$\cos 3\theta=4\cos^3\theta-3\cos\theta.$
Writing $z=\cos \theta$ we give it the form
\begin{equation}
\label{eq2-z}
z^3-\frac{3}{4}z-\frac{1}{4}\cos 3\theta=0,
\end{equation}
which clearly has the 
solutions 
\begin{equation}
\label{soleq2-zj}
z_j = \cos\left(\theta-\frac{2\pi j}{3}\right),\quad j=0,1,2.
\end{equation}
We see that equation \eqref{eq2-uz} will
be identical with \eqref{eq2-z} if $\displaystyle m=\sqrt{-\frac{4r}{3}}$ 
and \break
${\displaystyle \cos 3\theta=\frac{-27q}{\sqrt{-27r^3}}}$. Returning 
to the variable $h$ we see that 
\eqref{soleq2-zj} will lead to the 
solutions \eqref{eq2-25} to the original equation \eqref{eq2-22} and 
that we can choose $\theta$ as in 
\eqref{eq2-theta}.

Since $\theta\in (0,\frac{\pi}{3})$, we have
\[
-1<\cos\bigg(\theta-\frac{4\pi}{3}\bigg) < -\tfrac12 <
\cos\bigg(\theta-\frac{2\pi}{3}\bigg)<\tfrac12<\cos\theta <1.
\]
Hence $h_0>h_1>h_2$. Finally, note that with the notation from 
\eqref{eq2-thetaprim} we
have $\theta'=\frac{\pi}{3}-\theta$ and therefore $h'_0=h_2$,
$h«_1=h'_1$, $h'_2=h_0$ and thus $h'_{0}<h'_{1}<h'_{2}.$

\section{Solution of the main problem}
\label{sec3}
\setcounter{equation}{0}

\begin{theorem}
\label{thm3-1}
The set of equations \eqref{eq2-1} has exactly {\rm 20} solutions in
$\bC^{\,3}$. The first two solutions are the ``$\varepsilon$-solutions'':
\begin{equation}
\label{eq3-1}
c_0=c_1=c_2=\frac{2-p\pm\sqrt{p(p-4)}}{2}.
\end{equation}
The remaining {\rm 18} solutions can be obtained from the three solutions listed
below by the six transformations
\begin{eqnarray*}
(c_0,c_1,c_2) &\to & (c_k,c_{k+1},c_{k+2})\\
(c_0,c_1,c_2) &\to & 
\bigg(\frac{1}{c_k},\frac{1}{c_{k+1}},\frac{1}{c_{k+2}}\bigg)
\end{eqnarray*}
where $k=0,1,2$ and indices are computed modulo {\rm 3}. Put $u=\sqrt{p}$,
$v=\sqrt{p+4A+16}$ and 
$\theta={1\over 3}\Arccos\big(\frac{A}{2\sqrt{p}}\big)$. The three
solutions are $c^{(i})=(c_0^{(i)},c_1^{(i)},c_2^{(i)})$, $i=1,2,3,$ where
\begin{equation}
\label{eq3-2}
c_j^{(i)} =
\alpha^{(i)}+\beta^{(i)}\cos(\theta-\frac{2\pi}{3}j)
+\gamma^{(i)}\sin(\theta-\frac{2\pi}{3}j) 
\end{equation}
and

{\large

\renewcommand{\arraystretch}{1.5}
\begin{equation}
\label{eq3-3}
\left\{\begin{array}{lcl} \alpha^{(1)} &=& \frac12
\frac{pA-2p-2A}{p^2-3p-A}+i\frac{3\sqrt{3}}{2} \frac{\sqrt{p} 
\sqrt{p-4}
B}{p^2-3p-A} \\
\beta^{(1)} &=& -\frac12
\frac{\sqrt{p}(p-4)(A+2)}{p^2-3p-A}-i\frac{3\sqrt{3}}{2} 
\frac{\sqrt{p-4}
(p-2)B}{p^2-3p-A} \\
\gamma^{(1)} &=& -\frac{3\sqrt{3}}{2}
\frac{\sqrt{p}(p-4)B}{p^2-3p-A}+\frac{i}{2}\frac{\sqrt{p-4}(pA-2p-2A)}{p^2-3p-A}\,,

\end{array}\right.
\end{equation}

\begin{equation}
\label{eq3-4}
\left\{\begin{array}{lcl} \alpha^{(2)} &=& -\frac12
\frac{u^2-uv-4}{u^2+uv+2}+\frac{i}{2} \frac{u\sqrt{4+u-v} 
\sqrt{4-u+v}}{u^2+uv+2} \\
\beta^{(2)} &=& 
\frac{(A+2)u}{u^2+uv+2}+\frac{i}{4} \frac{(u^2+uv+4)\sqrt{4+u-v} 
\sqrt{4-u+v}}{u^2+uv+2} \\
\gamma^{(2)} &=& \frac{3\sqrt{3}Bu}{u^2+uv+2}+\frac{i}{4} 
\frac{(u^2-uv-4)\sqrt{u+v+4}\sqrt{u+v-4}}{u^2+uv+2}\,, 
\end{array}\right. 
\end{equation}

\begin{equation}
\label{eq3-5}
\left\{\begin{array}{lcl} \alpha^{(3)} &=& -\frac12
\frac{u^2+uv-4}{u^2-uv+2}-\frac{u}{2} \frac{\sqrt{u+v+4}
\sqrt{u+v-4}}{u^2-uv+2} \\ 
\beta^{(3)} &=& 
\frac{(A+2)u}{u^2-uv+2}-\frac14
\frac{(u^2-uv+4)\sqrt{u+v+4}\sqrt{u+v-4}}{u^2-uv+2} \\
\gamma^{(3)} &=& \frac{3\sqrt{3}Bu}{u^2-uv+2} + \frac14
\frac{(u^2+uv-4)\sqrt{4+u-v}\sqrt{4-u+v}}{u^2-uv+2}\,. 
\end{array}\right.
\end{equation}
\renewcommand{\arraystretch}{1.0}}
The solutions \eqref{eq3-2} given by \eqref{eq3-3} and \eqref{eq3-4} are
unimodular while the $\varepsilon$-solutions and the solution 
\eqref{eq3-2} given by
\eqref{eq3-5} are real. Hence of the {\rm 20} solutions {\rm 12} are unimodular 
and
{\rm 8} are real.

\end{theorem}

\begin{remark}{\rm \label{rem3-2} \ \\
(a) Of course the choice of a ``canonical'' solution among 
six possible ones is arbitrary. Our choice is motivated by a wish to 
give the asymptotic results in Section \ref{sec-asymptotics} a simple 
form.

\noindent(b) It follows from the proof of Theorem \ref{thm3-1} that the
transformation

\noindent $(c_0,c_1,c_2)\to
\left(\frac{1}{c_0},\frac{1}{c_1},\frac{1}{c_2}\right)$ can be 
obtained
just by changing the sign of the second term in the above formulas for
$\alpha^{(i)}$, $\beta^{(i)}$, and $\gamma^{(i)}$ .

\noindent(c) Since $u=\sqrt{p}$ and $v=\sqrt{p+4A+16}$ and
$|A|<2\sqrt{p}$, we have
\[
|u-4|<v<u+4
\]
which means that the numbers $u,v,4$ can be the lengths of the three
sides in a non-degenerate triangle. Hence the 4 square roots
\[
\sqrt{u+v+4},\qquad \sqrt{u+v-4},\qquad 
\sqrt{4+u-v},\qquad\sqrt{4-u+v}
\]
are well defined and strictly positive. Note also that
\begin{equation}
\label{Atouv}
A=\frac{v^2-u^2-16}{4}
\end{equation}
and
\begin{equation}
\label{Btouv}
B=\frac{1}{3\sqrt{3}} \sqrt{4p-A^2} =
\frac{\sqrt{u+v+4}\sqrt{u+v-4}\sqrt{4+u-v}\sqrt{4-u+v}}{12\sqrt{3}}.
\end{equation}}
\end{remark}
The proof of Theorem \ref{thm3-1} relies on Proposition \ref{prop2-1}
and the following 3 lemmas:

\begin{lemma}
\label{lemma3-2}
Let $a_0,a_1,a_2\in\bC$ and let $\theta\in\bR$. Then there are unique
numbers $\rho,\sigma,\tau\in\bC$ such that
\[
a_j=\rho+\sigma\cos\left(\theta-\frac{2\pi}{3}j\right)+\tau
\sin\left(\theta-\frac{2\pi}{3} j\right),\quad j=0,1,2.
\]
\end{lemma}   

\bproof
By an elementary computation one finds
\[
\det\left(\begin{array}{ccc} 1 & \cos\theta & \sin\theta\\
1 & \cos(\theta-\frac{2\pi}{3})&\sin(\theta-\frac{2\pi}{3})\\
1&\cos(\theta-\frac{4\pi}{3})&\sin(\theta-\frac{4\pi}{3})\end{array}\right)=-\frac{3\sqrt{3}}{2}.
\]
In particular the determinant is non-zero, which proves Lemma 
\ref{lemma3-2}.

\begin{lemma}
\label{lemma3-3}
Let $\theta\in\bR$ and let
$\alpha_1,\beta_1,\gamma_1,\alpha_2,\beta_2,\gamma_2\in\bC$, and put
\begin{eqnarray*}
c_j &=&
\frac{\alpha_1+\alpha_2}{2}+\frac{\beta_1+\beta_2}{2}\cos\left(\theta
-\frac{2\pi}{3}
j\right)+\frac{\gamma_1+\gamma_2}{2}\sin\left(\theta-\frac{2\pi}{3}
j\right)\\
\widetilde{c_j} &=&
\frac{\alpha_1-\alpha_2}{2}+\frac{\beta_1-\beta_2}{2}\cos\left(\theta
-\frac{2\pi}{3}
j\right)+\frac{\gamma_1-\gamma_2}{2}\sin\left(\theta-\frac{2\pi}{3}
j\right)
\end{eqnarray*}
for $j=0,1,2$. Then the following two conditions are equivalent
\medskip

{\rm(i)}   $c_0\tilde{c}_0 = c_1\tilde{c}_1 = c_2\tilde{c}_2=1$,
\medskip

{\rm(ii)}  $t_1=t_2=t_3=0$,

\medskip
\noindent where
\begin{eqnarray}
\label{eq3-6}
t_1 &=&
(\alpha_1^2-\alpha_2^2)+\tfrac12(\beta^2_1-\beta^2_2)+\tfrac12(\gamma_1^2-\gamma_2^2)-4,\\
\label{eq3-7}
t_2 &=&
2(\alpha_1\beta_1-\alpha_2\beta_2)+\tfrac12(\beta_1^2-\beta_2^2-\gamma_1^2+\gamma_2^2)\cos
3\theta +(\beta_1\gamma_1-\beta_2\gamma_2)\sin 3\theta,\\
\label{eq3-8}
t_3 &=&
2(\alpha_1\gamma_1-\alpha_2\gamma_2)+\tfrac12(\beta_1^2-\beta_2^2-\gamma_1^2+\gamma_2^2)\sin
3\theta -(\beta_1\gamma_1-\beta_2\gamma_2)\cos 3\theta.
\end{eqnarray}
\end{lemma}

\bproof
Put
\begin{eqnarray*}
f_j &=& c_j+\tilde{c}_j =
\alpha_1+\beta_1\cos\bigg(\theta-\frac{2\pi}{3}j\bigg)+\gamma_1\sin\bigg(\theta-\frac{2\pi}{3}j\bigg),\\
g_j &=& c_j-\tilde{c}_j =
\alpha_2+\beta_2\cos\bigg(\theta-\frac{2\pi}{3}j\bigg)+\gamma_2\sin\bigg(\theta-\frac{2\pi}{3}j\bigg).
\end{eqnarray*}
Then (i) is equivalent to
\[
f_j^2-g_j^2=4,\qquad j=0,1,2.
\]
By expressing $\cos^2\varphi$, $\sin^2\varphi$, 
$\cos\varphi\sin\varphi$
in terms of $\cos 2\varphi$, $\sin2\varphi$
($\varphi=\theta-\frac{2\pi}{3}j$) one finds
\begin{eqnarray*}
\label{eq3-9}
f^2_j &=&
\bigg(\alpha^2_1+\frac{\beta^2_1+\gamma^2_1}{2}\bigg)+2\alpha_1\beta_1\cos\bigg(\theta-\frac{2\pi}{3}
j\bigg)+2\alpha_1\gamma_1\sin\bigg(\theta-\frac{2\pi}{3}j\bigg)\\
&& \qquad
+\frac{\beta^2_1-\gamma^2_1}{2}\cos\bigg(2\theta-\frac{4\pi}{3}j\bigg)+
\beta_1\gamma_1\sin\bigg(2\theta-\frac{4\pi}{3}j\bigg).
\end{eqnarray*}
Using $\frac{4\pi}{3}j\equiv  -\frac{2\pi}{3}j\pmod{2\pi}$  one gets
\[
\left\{\begin{array}{lcl}
\cos\bigg(2\theta-\frac{4\pi}{3}j\bigg) &=&
\cos 3\theta\cos\bigg(\theta-\frac{2\pi}{3}j\bigg)+\sin 3\theta\sin\bigg(\theta-\frac{2\pi}{3}j\bigg)\\

\sin\bigg(2\theta-\frac{4\pi}{3}j\bigg)&=&
\sin 3\theta\cos\bigg(\theta-\frac{2\pi}{3}j\bigg)-\cos 3\theta\sin\bigg(\theta-\frac{2\pi}{3}j\bigg).
\end{array}\right.
\]
Hence
\begin{equation}
\label{f2coeffs}
f^2_j =
\rho_1+\sigma_1\cos\bigg(\theta-\frac{2\pi}{3}j\bigg)+\tau_1\sin\bigg(\theta-\frac{2\pi}{3}j\bigg),
\end{equation}
where
\renewcommand{\arraystretch}{1.5}
\[
\left\{\begin{array}{lcl}
\rho_1 &=& \alpha_1^2+\frac12(\beta^2_1+\gamma_1^2)\\
\sigma_1 &=&
2\alpha_1\beta_1+\frac{\beta^2_1-\gamma^2_1}{2}\cos 3\theta+\beta_1\gamma_1\sin 3\theta\\
\tau_1 &=&
2\alpha_1\gamma_1+\frac{\beta_1^2-\gamma_1^2}{2}\sin 3\theta-\beta_1\gamma_1\cos 3\theta.\end{array}\right.
\]
\renewcommand{\arraystretch}{1.0}
Similarly
\[
g_j^2 = \rho_2+\sigma_2
\cos\bigg(\theta-\frac{2\pi}{3}j\bigg)+\tau_2\sin\bigg(\theta-\frac{2\pi}{3}j\bigg),
\]
where
\[
\left\{\begin{array}{lcl}
\rho_2 &=& \alpha_2^2+\frac12(\beta^2_2+\gamma_2^2)\\
\sigma_2 &=&
2\alpha_2\beta_2+
2\alpha_2\beta_2\cos 3\theta)+\beta_2\gamma_2\sin 3\theta \\
\tau_2 &=&
2\alpha_2\gamma_2+
2\alpha_2\beta_2\sin 3\theta-\beta_2\gamma_2\cos 3\theta.
\end{array}\right.
\]
Since the coefficients in the decomposition
\[
f_j^2-g_j^2 =
(\rho_1-\rho_2)+(\sigma_1-\sigma_2)\cos\bigg(\theta-\frac{2\pi}{3}j\bigg)
+(\tau_1-\tau_2)\sin\bigg(\theta-\frac{2\pi}{3}j\bigg)
\]
are unique by Lemma \ref{lemma3-2}, we have $f_j^2-g_j^2=4$, 
$j=0,1,2$,
if and only if
\[
\rho_1-\rho_2=4,\quad \sigma_1-\sigma_2=0,\quad\mbox{and}\quad
\tau_1-\tau_2=0.
\]
This proves Lemma \ref{lemma3-3}.

\begin{lemma}
\label{lemma3-4}
Let $\theta\in\bR$ and let $c_0,c_1,c_2\in\bC\backslash\{0\}$. Put
\[
f_j=c_j+\frac{1}{c_j},\quad g_j=c_j-\frac{1}{c_j},\quad
h_j=\frac{c_{j+2}}{c_{j+1}}+\frac{c_{j+1}}{c_{j+2}},\quad
k_j=\frac{c_{j+2}}{c_{j+1}}-\frac{c_{j+1}}{c_{j+2}},
\]
where $j=0,1,2$
(counted modulo {\rm3}).

Let moreover
$\alpha_\nu,\beta_\nu,\gamma_\nu,\xi_\nu,\eta_\nu,\zeta_\nu$ 
($\nu=1,2$)
be the coefficients in the decompositions
\begin{equation}
\label{eq3-10}
\left\{\begin{array}{lcl} f_j &=&
\alpha_1+\beta_1\cos(\theta-\frac{2\pi}{3}j)+\gamma_1\sin(\theta-\frac{2\pi}{3}j)\\
g_j &=& 
\alpha_2+\beta_2\cos(\theta-\frac{2\pi}{3}j)+\gamma_2\sin(\theta-\frac{2\pi}{3}j)\end{array}\right.
\end{equation}

\begin{equation}
\label{eq3-11}
\left\{\begin{array}{lcl} h_j &=&
\xi_1+\eta_1\cos(\theta-\frac{2\pi}{3}j)+\zeta_1\sin(\theta-\frac{2\pi}{3}j)\\
k_j &=& 
\xi_2+\eta_2\cos(\theta-\frac{2\pi}{3}j)+\zeta_2\sin(\theta-\frac{2\pi}{3}j)\end{array}\right.
\end{equation}

Then
\begin{equation}
\label{eq3-12}
\left\{\begin{array}{lcl}
\xi_1 &=& \frac34(\alpha_1^2-\alpha_2^2)-1\\
\eta_1 &=& -\frac32(\alpha_1\beta_1-\alpha_2\beta_2)\\
\zeta_1 &=&
-\frac32(\alpha_1\gamma_1-\alpha_2\gamma_2)\end{array}\right.
\end{equation}
and
\begin{equation}
\label{eq3-13}
\left\{\begin{array}{lcl}
\xi_2 &=& \frac{\sqrt{3}}{4}(\beta_2\gamma_1-\beta_1\gamma_2)\\
\eta_2 &=& \frac{\sqrt{3}}{2}(\gamma_2\alpha_1-\gamma_1\alpha_2)\\
\zeta_2 &=&
\frac{\sqrt{3}}{2}(\alpha_2\beta_1-\alpha_1\beta_2).\end{array}\right.
\end{equation}
\end{lemma}

\bproof
Clearly
\[
c_j = \tfrac12(f_j+g_j), \qquad \tfrac{1}{c_j}=\tfrac12(f_j-g_j).
\]
Hence
\begin{eqnarray*}
h_j &=& \tfrac12(f_{j+1}f_{j+2}-g_{j+1}g_{j+2})\\
k_j &=& \tfrac12(f_{j+1}g_{j+2}-g_{j+1}f_{j+2}).
\end{eqnarray*}
By expressing $\cos(\theta-\frac{2\pi}{3})$,
$\sin(\theta-\frac{2\pi}{3})$, $\cos(\theta-\frac{2\pi}{3})$, and
$\sin(\theta-\frac{4\pi}{3})$ as linear combinations of $\cos\theta$
and $\sin\theta$ one gets
\begin{equation}
f_1f_2 =
(\alpha^2_1-\frac{\beta^2_1+\gamma^2_1}{4})-\alpha_1\beta_1\cos\theta-\alpha_1\gamma_1\sin\theta 
+ \frac{\beta^2_1-\gamma^2_1}{2}\cos
2\theta+\beta_1\gamma_1\sin 2\theta. 
\end{equation}
Using now \eqref{f2coeffs} from the proof of Lemma \ref{lemma3-3}, we 
have
\[
f_1f_2-f_0^2 = -\tfrac34 
(\beta^2_1+\gamma_1^2)-3\alpha_1\beta_1\cos\theta-3\alpha_1\gamma_1\sin\theta.
\]
Repeating the same argument with $\theta-\frac{2\pi}{3}j$ instead
of $\theta$, we have
\begin{equation}
\label{eq3-14}
f_{j+1}f_{j+2}-f_j^2=-\tfrac34(\beta_1^2+\gamma_1^2)
-3\alpha_1\beta_1\cos\bigg(\theta-\frac{2\pi}{3}j\bigg)
-3\alpha_1\gamma_1\sin\bigg(\theta-\frac{2\pi}{3}j\bigg)
\end{equation}
and in the same way we have
\begin{equation}
\label{eq3-15}
g_{j+1}g_{j+2}-g_j^2=-\tfrac34(\beta_2^2+\gamma_2^2)
-3\alpha_2\beta_2\cos\bigg(\theta-\frac{2\pi}{3}j\bigg)-3\alpha_2\gamma_2\sin\bigg(\theta
-\frac{2\pi}{3}j\bigg).
\end{equation}
By the definition of $f_j$ and $g_j$ we have
\begin{equation}
\label{eq3-16}
f_j^2-g_j^2=\bigg(c_j+\frac{1}{c_j}\bigg)^2-\bigg(c_j-\frac{1}{c_j}\bigg)^2=4.
\end{equation}
Hence, by \eqref{eq3-14}, \eqref{eq3-15}, and \eqref{eq3-16}
\begin{eqnarray*}
2h_j &=& f_{j+1}f_{j+2}-g_{j+1}g_{j+2}\\
&=&
4-\tfrac34(\beta_1^2+\gamma_1^2-\beta_2^2-\gamma_2^2)-3(\alpha_1\beta_1
-\alpha_2\beta_2)\cos(\theta-\tfrac{2\pi}{3}j)\\
&\quad\quad&-
3(\alpha_1\gamma_1
-\alpha_2\gamma_2)\sin(\theta-\tfrac{2\pi}{3}j).
\end{eqnarray*}
By uniqueness of this decomposition (Lemma \ref{lemma3-2}) we can read
off the coefficients $\xi_1,\eta_1,\zeta_1$ in \eqref{eq3-11} namely
\begin{eqnarray*}
\xi_1 &=& 2-\tfrac38(\beta_1^2+\gamma_1^2-\beta_2^2-\gamma_2^2),\\
\eta_1 &=& -\tfrac32(\alpha_1\beta_1-\alpha_2\beta_2),\\
\zeta_1 &=& -\tfrac32(\alpha_1\gamma_1-\alpha_2\gamma_2).
\end{eqnarray*}
However by \eqref{eq3-6} in Lemma \ref{lemma3-3}, we have
\[
(\alpha_1^2-\alpha_2^2)+\tfrac12(\beta_1^2-\beta_2^2)+\tfrac12(\gamma^2_1-\gamma^2_2)=4.
\]
Hence the above formula for $\xi_1$ can be changed to
\[
\xi_1 = \tfrac34(\alpha_1^2-\alpha_2^2)-1.
\]
This proves \eqref{eq3-12}. A similar but much simpler computation 
gives
\begin{eqnarray*}
k_j &=& \tfrac12(f_{j+1}g_{j+2}-f_{j+2}g_{j+1})\\
&=&
\frac{\sqrt{3}}{4}(\beta_2\gamma_1-\beta_1\gamma_2)+\frac{\sqrt{3}}{2}(\gamma_2\alpha_1
-\gamma_1\alpha_2)\cos\bigg(\theta-\frac{2\pi}{3}j\bigg)\\
&&\qquad 
+\frac{\sqrt{3}}{2}(\alpha_2\beta_1-\alpha_1\beta_2)\sin\bigg(\theta-\frac{2\pi}{3}j\bigg),
\end{eqnarray*}
which proves \eqref{eq3-13}.

\noindent{\bf Proof of Theorem \ref{thm3-1}:} 
Assume that $(c_0,c_1,c_2)$ is a solution to the set of equations
\eqref{eq2-1}. By Proposition \ref{prop2-1}, the numbers
\[
h_j = \frac{c_{j+2}}{c_{j+1}}+\frac{c_{j+1}}{c_{j+2}},\qquad j=0,1,2,
\]
must be of the form
\begin{equation}
\label{eq3-17}
h_j = \xi_1+\eta_1\cos\bigg(\theta-\frac{2\pi}{3}j\bigg),\qquad 
j=0,1,2,
\end{equation}
where $(\xi_1,\eta_1)$ is one of the four pairs
$(\xi_1^{(i)},\eta_1^{(i)})$, $i=0,1,2,3$, listed in
\eqref{eq2-4}--\eqref{eq2-7}. For $i=0$, we have $\xi_1=2$ and
$\eta_1=0$. Hence $h_0=h_1=h_2=2$ which implies that $c_0=c_1=c_2$, 
and
in this case the only solutions to \eqref{eq2-1} are the 2
``$\varepsilon$-solutions'' from \cite{pisa}, namely
\[
c_0=c_1=c_2=\frac{2-p\pm\sqrt{p(p-4)}}{2}.
\]
For $i=1,2,3$ we can compute the numbers $c_j$ from $(\xi_1,\eta_1)$ 
by
Lemma \ref{lemma3-4}. Define
\[
f_j=c_j+\frac{1}{c_j},\quad g_j=c_j-\frac{1}{c_j},\quad
h_j=\frac{c_{j+2}}{c_{j+1}}+\frac{c_{j+1}}{c_{j+2}},\quad
k_j=\frac{c_{j+2}}{c_{j+1}}-\frac{c_{j+1}}{c_{j+2}}
\]
as in Lemma \ref{lemma3-4}, and let
$\alpha_\nu,\beta_\nu,\gamma_\nu,\xi_\nu,\eta_\nu,\zeta_\nu$,
$\nu=1,2$ be the coefficients in the decompositions \eqref{eq3-10}
and \eqref{eq3-11}. Note that by Lemma \ref{lemma3-2} this new definition of $\xi_1$ and $\eta_1$ 
is
consistent with \eqref{eq3-17}. Moreover $\zeta_1=0$ by 
\eqref{eq3-17}.

By \eqref{eq2-1}
\[
f_j=-\frac{p-4}{3}-\frac{p+A+1}{9} h_j-\frac{2p-A-9B-4}{18}
h_{j+1}-\frac{2p-A-9B-4}{18} h_{j+2}.
\]
Since
\begin{eqnarray*}
h_0 &=& \xi_1+\eta_1\cos\theta,\\
h_1 &=& \xi_1+\eta_1\bigg(-\tfrac12\cos\theta +
i\frac{\sqrt{3}}{2}\sin\theta\bigg),\\
h_2 &=& \xi_1+\eta_1\bigg(-\tfrac12\cos\theta 
-i\frac{\sqrt{3}}{2}\sin\theta\bigg), 
\end{eqnarray*}
we have
\[
f_0 =
\left(-\frac{p-4}{3}-\frac{p-1}{3}\xi_1\right)-\frac{A+2}{6}\eta_1\cos\theta-\frac{\sqrt{3}}{2}B\eta_1.
\]
Repeating the same computation with $\theta$ replaced by
$\eta-\frac{2\pi}{3}j$, we get that the coefficients
$\alpha_1,\beta_1,\gamma_1$ in the decomposition
\[
f_j=\alpha_1+\beta_1\cos\bigg(\theta-\frac{2\pi}{3}j\bigg)+\gamma_1\sin\bigg(\theta-\frac{2\pi}{3}j\bigg)
\]
are given by
\begin{equation}
\label{eq3-18}
\left\{ \begin{array}{lcl} \alpha_1 &=&
-\frac{p-4}{3}-\frac{p-1}{3}\xi_1\\
\beta_1 &=& -\frac{A+2}{6} \eta_1\\
\gamma_1 &=& -\frac{\sqrt{3}}{2}B. \end{array}\right.
\end{equation}
Provided $\alpha_1^2-\frac43(\xi_1+1)\ne 0$ we then get from
\eqref{eq3-12}
\begin{equation}
\label{eq3-19}
\left\{ \begin{array}{lcl} \alpha_2 &=&
\pm\sqrt{\alpha_1^2-\frac43(\xi_1+1)}\\
\beta_2 &=& \frac{1}{\alpha^2}(\alpha_1\beta_1+\frac23\eta_1)\\
\gamma_2 &=& \frac{1}{\alpha_2}(\alpha_1\gamma_1+\frac23\zeta_1).
\end{array}\right.
\end{equation}
Inserting the values $(\xi_1^{(i)},\eta_1^{(i)})$, $i=1,2,3$ from
\eqref{eq2-20} in \eqref{eq3-18} we find that
$\alpha_1^2-\frac43(\xi_1+1)\ne 0$ in all the cases $i=1,2,3$. Hence 
the
numbers $\alpha_1$, $\beta_1$, $\gamma_1$, $\alpha_2$, $\beta_2$,
$\gamma_2$ given by \eqref{eq3-18} and \eqref{eq3-19} are unique up to
simultaneous sign change of $(\alpha_2,\beta_2,\gamma_2)$. For 
$i=1,2$,
\[
\alpha_2=\pm\sqrt{\alpha_1^2-\frac43(\xi_1+1)}
\]
is purely imaginary, and we choose the solution with
$\Im(\alpha_2^{(i)})>0$ ($i=1,2$). For $i=3$, $\alpha_2$ is real and 
we
choose the solution with $\sign(\alpha_2^{(3)})=-\sign(u^2-uv+2)$. It 
is
now easy to compute $\alpha_1$, $\beta_1$, $\gamma_1$, $\alpha_2$,
$\beta_2$, $\gamma_2$ explicitly from \eqref{eq2-20} in the 3 cases
$i=1,2,3$. One finds
{\large  \renewcommand{\arraystretch}{1.5}
\begin{equation}
\label{eq3-20}
\left\{ \begin{array}{lclclcl} \alpha_1^{(1)} &=&
\frac{pA-2p-2A}{p^2-3p-A},\quad & \alpha_2^{(1)} &=&
i\frac{3\sqrt{3}\sqrt{p}\sqrt{p-4}B}{p^2-3p-A}\\
\beta_1^{(1)} &=& -\frac{\sqrt{p}(p-4)(A+2)}{p^2-3p-A}, & 
\beta_2^{(1)}
&=& -i\frac{3\sqrt{3}\sqrt{p-4}(p-2)B}{p^2-3p-A}\\
\gamma_1^{(1)} &=& -\frac{3\sqrt{3}\sqrt{p}(p-4)B}{p^2-3p-A},\quad &
\beta_3^{(1)} &=& i\frac{\sqrt{p-4}(pA-2p-2A)}{p^2-3p-A}. 
\end{array}\right.
\end{equation}

\begin{equation}
\label{eq3-21}
\left\{ \begin{array}{lclclcl} \alpha_1^{(2)} &=&
-\frac{u^2-uv-4}{u^2+uv+2},\quad & \alpha_2^{(2)} &=&
i\frac{u\sqrt{4+u-v}\sqrt{4-u+v}}{u^2+uv+2}\\
\beta_1^{(2)} &=& \frac{2(A+2)u}{u^2+uv+2}, & \beta_2^{(2)}
&=& \frac{i}{2}\frac{(u^2+uv+4)\sqrt{4+u-v}\sqrt{4-u+v}}{u^2+uv+2}\\
\gamma_1^{(2)} &=& \frac{6\sqrt{3}Bu}{u^2+uv+2},\quad &
\gamma_2^{(2)} &=& 
\frac{i}{2}\frac{(u^2-uv+4)\sqrt{u+v+4}\sqrt{u+v-4}}{u^2+uv+2}. 
\end{array}\right.
\end{equation}

\begin{equation}
\label{eq3-22}
\left\{ \begin{array}{lclclcl} \alpha_1^{(3)} &=&
-\frac{u^2+uv-4}{u^2-uv+2},\quad & \alpha_2^{(3)} &=&
-\frac{u\sqrt{u+v+4}\sqrt{u+v-4}}{u^2-uv+2}\\
\beta_1^{(3)} &=& \frac{2(A+2)u}{u^2-uv+2}, & \beta_2^{(3)}
&=& -\frac12\frac{(u^2-uv+4)\sqrt{u+v+4}\sqrt{u+v-4}}{u^2-uv+2}\\
\gamma_1^{(3)} &=& \frac{3\sqrt{3}Bu}{u^2-uv+2},\quad &
\gamma_2^{(3)} &=& 
\frac12\frac{(u^2+uv-4)\sqrt{4+u-v}\sqrt{4-u+v}}{u^2-uv+2}. 
\end{array}\right.
\end{equation}}
\renewcommand{\arraystretch}{1.0}
Since
\begin{equation}
\label{eq3-23}
c_j=\tfrac12(f_j+g_j)=\frac{\alpha_1+\alpha_2}{2}+\frac{\beta_1+\beta_2}{2}\cos\bigg(\theta-\frac{2\pi}{3}j\bigg)+\frac{\gamma_1+\gamma_2}{2}\sin\bigg(\theta-\frac{2\pi}{3}j\bigg),
\end{equation}
we obtain \eqref{eq3-2} with $\alpha^{(i)},\beta^{(i)},\gamma^{(i)}$
given by \eqref{eq3-3}, \eqref{eq3-4} and \eqref{eq3-5}.

We still have to check that the $(c_0^{(i)},c_1^{(i)},c_2^{(i)})$ 
given
by \eqref{eq3-2}--\eqref{eq3-5} actually are solutions to
\eqref{eq2-1}. From Lemma \ref{lemma3-3} and Lemma \ref{lemma3-4} it
follows that the only thing left to check is that $c_j\ne 0$, 
$j=0,1,2$
and that
\begin{equation}
\label{eq3-24}
\frac{1}{c_j}=\frac{\alpha_1-\alpha_2}{2}+
\frac{\beta_1-\beta_2}{2}\cos\bigg(\theta-\frac{2\pi}{3}j\bigg)
+\frac{\gamma_1-\gamma_2}{2}\sin\bigg(\theta-\frac{2\pi}{3}
j\bigg),
\end{equation}
which is equivalent to checking that the numbers $t_1,t_2,t_3$ listed in
\eqref{eq3-6}--\eqref{eq3-8} are zero.

Using
\[
\cos 3\theta=\frac{A}{2\sqrt{p}},\quad \sin 3\theta =
\frac{\sqrt{4p}-A^2}{2\sqrt{p}}=\frac{3\sqrt{3} B}{2\sqrt{p}}
\]
it is elementary to check by MAPLE or MATHEMATICA that $t_1=t_2=t_3=0$
in each of the 3 cases \eqref{eq3-20}, \eqref{eq3-21} and 
\eqref{eq3-22}
above. It is also possible to avoid a case by case check by relating
$t_1$, $t_2$ and $t_3$ to the polynomials $p_1,p_2,p_3,p_4$ used in the
proof of Proposition \ref{prop2-1} (see Remark \ref{rem3-5} below).

Finally we have to show that we have found 20 distinct solutions: 
Since
$\eta_1^{(i)}\ne 0$, $i=1,2,3$, the 3 solutions given by
\eqref{eq3-2}--\eqref{eq3-5} are distinct from the two
$\varepsilon$-solutions. This also implies that in each of the 3 
cases,
the 6 solutions given by
\begin{equation}
\label{eq3-25}
\left\{\begin{array}{ll}
(c_j,c_{j+1},c_{j+2})\qquad & j=0,1,2\\
(\frac{1}{c_j},\frac{1}{c_{j+1}},\frac{1}{c_{j+2}}, & 
j=0,1,2\end{array}\right.
\end{equation}
are all distinct. To check that there is no overlap between these 3
groups of 6 solutions it is sufficient to check that the 3 numbers
$s_1^{(i)}=3\xi_1^{(i)}$ are distinct because
\[
s_1=h_0+h_1+h_2=\frac{c_2}{c_1}+\frac{c_0}{c_2}+\frac{c_1}{c_0}+\frac{c_1}{c_2}+\frac{c_2}{c_0}+\frac{c_0}{c_1}
\]
is invariant under the 6 transformations listed in \eqref{eq3-25}. 
From
\eqref{eq2-19}
\begin{eqnarray*}
s_1^{(1)} &=& \frac{18p-3p^2-6A}{p^2-3p-A}\\
\left.\begin{array}{ll} s_1^{(2)}\\ s_1^{(3)} \end{array}\right\} &=&
\frac{-6pA-27-12\pm 9\sqrt{p(p+4A+16)}}{2(pA+2p-1)}.
\end{eqnarray*}
Clearly $s_1^{(2)}\ne s_1^{(3)}$, since $p+4A+16>0$ by Remark
\ref{rem2-2}. Moreover $s_1^{(2)}$ and $s_1^{(3)}$ are the two zeros 
of
the polynomial $r$ from \eqref{eq2-18}:
\[
r(s_1)=(pA+3p-1)s_1^2+(6pA+27p+12)s_1 +(9pA+54p-36).
\]
We get
\[
r(s_1^{(1)}) = 81\frac{(2p-A-4)(4p-A^2)}{(p^2-3p-A)^2},
\]
but $2p-A-4 > 2p-2\sqrt{p}-4=2(\sqrt{p}+1)(\sqrt{p}-2)>0$, and
$4p-A^2=27B^2>0$. Hence $s_1^{(1)}\ne s_1^{(2)}$ and $s_1^{(1)}\ne
s_1^{(3)}$. Therefore we have found altogether $2+3\cdot 6=20$
solutions. By \eqref{eq3-24}, passing from $c_j^{(i)}$ to
$\frac{1}{c_j^{(i)}}$ in \eqref{eq3-2} corresponds to a change of sign
of $\alpha_2,\beta_2$ and $\gamma_2$. Hence the 12 solutions generated
by \eqref{eq3-3}, \eqref{eq3-4} and the transformations \eqref{eq3-25}
are all unimodular while the remaining 8 solutions clearly are real.

This completes the proof of Theorem \ref{thm3-1}.


\begin{remark}
\label{rem3-5}
{\rm We sketch here a different proof of $t_1=t_2=t_3=0$ for the values of
$\alpha_1,\beta_1,\gamma_1,\alpha_2,\beta_2,\gamma_2$
listed in
\eqref{eq3-20}--\eqref{eq3-22}:

By \eqref{eq3-18} and \eqref{eq3-19},
$\alpha_1,\beta_1,\gamma_1,\alpha_2,\beta_2,\gamma_2$ can be expressed
in terms of $(\xi_1,\eta_1)$ and hence $t_1,t_2,t_3$ given by
\eqref{eq3-6}--\eqref{eq3-7} can be expressed in terms of 
$\xi_1,\eta_1$, 
and $\theta$. Next we observe that if
\[
h_j=\xi_1+\eta_1\cos\bigg(\theta-\frac{2\pi}{3}j\bigg),\quad j=0,1,2,
\]
then
\begin{eqnarray*}
s_1 &=& h_0+h_1+h_2=3\xi_1,\\
s_2 &=& h_0h_1+h_1h_2+h_2h_0=3\xi_1^2-\frac34 \eta_1^2,\\
s_3 &=& h_0h_1h_2 = \xi_1^3-\tfrac34\xi_1\eta_1+\tfrac14\eta_1^3\cos
3\theta,\\
a &=& (h_0-h_1)(h_1-h_2)(h_2-h_0)=-\frac{3\sqrt{3}}{4} \eta_1\sin 
3\theta.
\end{eqnarray*}
Inserting this into the 4 polynomials $p_i=p_i(s_1,s_2,s_3,a)$ from 
the
proof of Proposition \ref{prop2-1} and comparing these new formulas 
for
$p_1,p_2,p_3$ and $p_4$ with the formulas found above for
$t_1,t_2,t_3$ one discovers after some work that
\renewcommand{\arraystretch}{1.5}
\begin{eqnarray*}
 t_1 &=& \frac{4(p_3-p_4)}{27\alpha_2^2}\\
 t_2 &=&
 -\frac{4(3p_1+(\xi_1-p\xi_1-p+4)p_2+(2\xi_1-1)p_3+(\xi_1+4)p_4}{27\xi_1\alpha_2^2}\\
 t_3 &=&
 \frac{4\Big(\xi_1+p\xi_1+A\xi_1+p+A-2\Big)\Big((p-1)p_1-2p_2-p_3-2p_4\Big)}{3\sqrt{3} 
B\xi_1\alpha_2^2}
\end{eqnarray*}
\renewcommand{\arraystretch}{1.0}
and since $(\xi_1^{(i)},\eta_1^{(i)})$, $i=1,2,3$ were found by 
solving
the equations $p_1=p_2=p_3=p_4=0$, it follows that $t_1=t_2=t_3=0$ in all
three cases.}
\end{remark}

\section{Corollaries of the main result (Leaving the simple case)}
\label{sec-results}
\setcounter{equation}{0}
In this section we will formulate and prove various consequences of 
the main result; in particular we will identify all bi-unimodular 
$p$-sequences and cyclic $p$-roots of index 3.
We will give 
the $c^{(i)}$ names:
\newtheorem {canonic}{Definition} [section]
\begin{canonic}
\label{canonic}
We denote as the {\rm first, second and third canonical solution} the
solutions $c^{(1)}$, $c^{(2)}$, and $c^{(3)}$ defined in Theorem {\rm\ref{thm3-1}}.
\end{canonic}

We will start by presenting all bi-unimodular $p$-sequence of index 
$3$ (cf. Definition \ref{general}). Recall that $\k=\exp ({2\pi i\over p}).$ 
\newtheorem{allbiuni}[canonic]{Proposition}
\begin {allbiuni}
  Let $p$ be a prime $\equiv 1 \pmod{6}$, and let $x$ be a
  bi-unimodular $p$-sequence of index $3$. Then there are a complex
  number $b$ of modulus one and integers $~r$ and $l$ such that $x$ is
  given by $x_l=b$ and $x_j=b\cdot\k^{rj}\cdot c_k~ {\rm when}~0\ne j-l\in
  G_k~~~(k=0,1,2)$, where $c=(c_{0},c_{1},c_{2})$ is one of the {\rm 12} 
  solutions to \eqref{eq2-1} coming from the the first or second canonical 
  solution $ c^{(1)},  c^{(2)}$,
  as described in Theorem {\rm\ref{thm3-1}}. If
  $p\ne 7$, there are $12p^2$ different {\rm normalized} 
  bi-unimodular $p$-sequences of index $3$ (i.e. with $x_0=1$).
  There are $336$ different normalized bi-unimodular $7$-sequences.
  Of these, $6\cdot 7^2$ come from the second canonical solution, whereas 
  only $6\cdot 7$ come from the first canonical solution.
  The last-mentioned sequences can be uniquely written in the form 
  $\displaystyle 
  x_j=\k^{m\cdot j^2 +n j}$,
  where $m$ and $n\in\bZ_{7}$ and $m\ne0$.
\label{allbiuni}
\end{allbiuni}

Next we formulate our result as a theorem bearing on cyclic $p$-roots
rather than on bi-unimodular $p$-sequences:

\newtheorem{allcyk}[canonic]{Proposition}
\begin{allcyk}
  Let $p$ be a prime $\equiv 1 \pmod{6}$, and let $z=\row z$ be a cyclic
  $p$-root of index $3$. Then  there are integers $r$ and $l$ such that $z$ is given by
  $z_j=\k^r\cdot c_{k}/c_\kappa$ when $j+1-l\in
  G_k$ and $j-l\in G_{\kappa}$, where $c=(c_{0},c_{1},c_{2})$ is one of 
  the {\rm 20} solutions to \eqref{eq2-1} as described in Theorem {\rm\ref{thm3-1}}.
  If $p\ne 7$, there are $20p^2$ different cyclic
  $p$-roots of index $3$, ($2p^2$ of which being in fact of index {\rm 1}). 
  There are only {\rm 434} different cyclic {\rm 7}-roots of index 
  {\rm 3}. Of 
  these, {\rm 42} come from the first canonical solution. These ``Gaussian'' 
  cyclic {\rm 7}-roots can be uniquely written in the form 
  $z_j=\k^{mj+n}$ 
   where $m$ and $n\in\bZ_{7}$ and $m\ne0$.
  \label{allcyk}
  \end{allcyk}

\medskip
\noindent{\bf Proof of Proposition \ref{allbiuni} and  Proposition \ref{allcyk}} 
The first statements in these theorems are obvious reformulations of 
Theorem \ref{thm3-1}
in terms of the concepts introduced in Section \ref{sec-notation}, and we 
leave it to the reader to check this. We will only prove the statements 
about the number of different normalized bi-unimodular 
sequences of index 3 (NBUS3), the number of different cyclic 
$p$-roots of index 3, and the explicit forms given in the first canonical 
case for $p=7$. 

We start with the last topic. Since the 42 possible $\k$-exponents 
in the $z_j$-formula in Proposition \ref{allcyk} form the set of 
all differences (as functions of $j$) of those  
in the $x_j$-formula in Proposition \ref{allbiuni}, it suffices to 
consider the latter (cf. \eqref{eq:xtoz} and Proposition 
\ref{xiffz}). We start by  taking $m=1$ and $n=0$, which gives
$x=(1,\k,\k^4,\k^2,\k^2,\k^4,\k).$
Since for $p=7$ we have $G_{0}=\{1,6\}$, $G_{1}=\{3,4\}$, and 
$G_{2}=\{2,5\}$, this means that this particular $x$ is in fact  
simple of index 3 with $c_{0}=\k$, 
$c_{1}=\k^2$, and $c_{2}=\k^4$ (cf. Definition \ref{simple}). We 
claim that this $c=(c_{0},c_{1},c_{2})$ is one of the six solutions 
coming from $c^{(1)}$ in Theorem \ref{thm3-1}. To prove this, we 
calculate $h_{0}=\frac{c_{1}}{c_{2}}+\frac{c_{2}}{c_{1}}=\k^2+\k^{-2}$, 
$h_{1}=\frac{c_{2}}{c_{0}}+\frac{c_{0}}{c_{2}}=\k^3+\k^{-3}$,
and $h_{2}=c_{0}/c_{1}+c_{1}/c_{0}=\k+\k^{-1}$. Thus, using the 
relation $1+\k+\k^2+\k^3+\k^4+\k^5+\k^6=0$, we get 
$s_{1}=h_{0}+h_{1}+h_{2}=-1$, $s_{2}=h_{0}h_{1}+ 
h_{1}h_{2}+h_{2}h_{0}=-2$, $s_{3}=h_{0}h_{1}h_{2}=1$, 
and                                
$a=(h_{1}-h_{0})(h_{2}-h_{1})(h_{0}-h_{02})=-7$. Since these values 
agree with those of $s_{1}^{(1)},s_{2}^{(1)},s_{3}^{(1)},$ and $a^{(1)}$ 
in \eqref{sa7}, our last claim is proved. 

Next we 
keep $n=0$ but consider a general $m$. But all we have used about $\k$ in our 
calculations is that $\k$ is a primitive seventh root of unity. So is 
$\k^m$. Thus, 
 $\displaystyle x_{j}=\k^{m\cdot j^2}=(\k^m)^{j^2}$ will also give a 
simple  bi-unimodular 7-sequence of index 3. Of course the six possibilities for $m$ 
correspond to the six transformations mentioned in Theorem 
\ref{thm3-1}.
Finally, taking a
general $n$, we see by Definition \ref{general} (with $l=0$ and $h=n$ ) 
that all our $x$ are bi-unimodular 7-sequences of index 3. Clearly 
they are normalized.

It is clear that the 42 normalized bi-unimodular 7-sequences of index 
3 we have found are different. Next we show that no other normalized 
bi-unimodular 7-sequence comes from the first canonical case. All we 
have to prove is that taking $l\ne 0$ in Definition \ref{general} does 
not give anything new when $y$ is a simple bi-uninormal sequence of index 3 given by 
$y_{k}=\k^{mj^2}$. But this is
trivial, since Definition \ref{general} gives the unnormalized bi-uninormal sequence $x$ of index 3 
 defined by $x_{j}=\k^{hj+m(j-l)^2} 
=\k^{ml^2+mj^2-2mlj}$ which is normalized through division by 
$x_{0}=\k^{ml^2}$ and becomes $\k^{mj^2-2mlj}=\k^{(-2ml)j}y_{j}$, which is of the 
desired form. 

It remains to prove that the numbers of different NBUS3:s and 
different cyclic $p$-roots of index 3 given in our two propositions are 
correct, that is that no such ``collapse'' occurs except in the 
first canonical case for for $p=7$. Recall that in the end of the proof of
Theorem \ref{thm3-1} we showed that 
all the 20 solutions to the main problem are different. We now have to 
extend this from the simple to the general case and we start by considering 
the $\varepsilon$-solutions. Every corresponding NBUS3 $x$ has the form 
$x_{j}=d_{j}\k^{rj}$ with $r\in\bZ$ and 
$d=(1,\varepsilon,\varepsilon,\ldots,\varepsilon,\varepsilon)$ or
$d=(\ldots,1,1,\varepsilon,1,1,\ldots)$ with 
$\varepsilon=\big(2-p\pm\sqrt{p(p-4)}\,\,\big)/2.$
These $p^2$ NBUS3:s are clearly distinct. 

Let us when $r\ne0$ and $l$ are in $\bZ_p$ and 
$c=(c_{0},c_{1},c_{2})\in \bC^{\,3}$ is one of the 20  
solutions mentioned in 
\ref{thm3-1}, 
define $x(r,l,c)$ as the NBUS3  $x=(x_0, x',\ldots, 
x_{p-1})$ given by the formulas:
\begin {equation}
x_j=b\k^{rj}(c)_k \,\,{\rm when}\,\, 0\ne j-l\in G_k,
\label{jnotl}
\end{equation}
\begin {equation}
 x_l=b\k^{rl},
\label{jisl}
\end{equation}
where $b$ is determined by the normalization 
\begin {equation}
x_0=1.
\label{bnorm}
\end{equation}
 
Let us consider two coinciding NBUS3:s,  $x(r',l',c')=x(r'',l'',c'')$
which do not satisfy all the three equalities $r'=r'',l'=l'',c'=c''$. 
We denote the two $b$:s defined by 
(\ref{bnorm}) by $b'$ and $b''$, respectively.
We start by considering the possibility that
$l''=l'$. Denote the common value by $l$ and fix a $k$. 
From
(\ref{jnotl}) follows that
$b'\k^{r' j}(c')_k=b''\k^{r'' j}(c'')_k$ 
if $j-l\in G_k$ and thus for
at least two different non-zero $j$, which leads to 
$r'=r''.$ Then (\ref{jisl})  gives $b'=b''$. Now \eqref{jnotl} implies 
that we have also $c'=c'',$ which is
against our hypothesis that at least one of
$r,l$ and $c$ differs between the two NBUS3:s.
 
Thus we have
$l'\ne l''$.  Let us now  suppose that $r'=r''$ (and $l'\ne l''$). 
Denote the common
$r$-value by $r$. Choose $j_1$ such that $j_{1}\ne l'$ and $j_{1}\ne l''$ and 
define $k_1$ and $k_2$ by
\begin{equation} 
(j_1-l')\in G_{k_1},\,\,\,\,(j_1-l'')\in G_{k_2}.
\label{k1k2def}
\end{equation}
Consider the set $F:=\{j-l'';(j-l')\in G_{k_1}\}\cap G_{k_2}$. Taking 
$d=l''-l'$ in (\ref{eq:transnos}), we see that if $d\in G_a$, 
then the cardinality of $F$ is a transition number:  
$\sharp(F)=n_{k_1-a,k_2-a}$. By (\ref{nlist}), all transition numbers are 
$\le s-1$, and since $\sharp(G_{k_1})=s$, there is at least one
 $j_2$ and one $k_3\ne k_2$ such that
\begin{equation} 
(j_2-l')\in G_{k_1}\,\,\,\,{\rm and}\,\,\,\,(j_2-l'')\in G_{k_3}
\label{k1k2diff}
\end{equation}
Now from (\ref{k1k2def}) and (\ref{k1k2diff}) 
follows that (\ref{jnotl}) with $j=j_1$ and with $j=j_2$ gives
$$b'\k^{rj_1}c_{k_1}'=b''\k^{rj_1}c_{k_2}'',$$
$$b'\k^{rj_2}c_{k_1}'=b''\k^{rj_2}c_{k_3}''.$$
This leads to $c''_{k_{2}}=c''_{k_{3}.}$ 
Then it follows from from Remark \ref{anot0} that $c''$ is an $\varepsilon$-solution. 
Since $c'$ and $c''$ play the same part in our situation, the same 
must be true for $c'$. But we know already that there is no
internal collapse among the NBUS3:s coming from 
$\varepsilon$-solutions, so the case $r'=r''$ also leads to a contradiction.
Now we know that $r'\ne r''$ and $l'\ne l''.$
From (\ref{jisl}) and (\ref{jnotl}) with $j=l'$ we get
\begin{equation}
x_{l'}=b'\k^{r' l'} =b''\k^{r'' l'}c_k'',
\label{jisl1}
\end{equation}
where $k$ is determined by $(l''-l')\in G_k$.
Since $G_k$
has at least two elements we can choose $j\ne l''$ with $(j-l')\in 
G_k$. 
For this $j$ we get from (\ref{jnotl})
\begin{equation}
x_j=b'\k^{r' j}c_k'=b''\k^{r'' j}c_k''.
\label{jnotl1}
\end{equation}
From (\ref{jisl1}) and (\ref{jnotl1}) we get by division
$$c_k'=\k^{(r''-r')(j-l')}.$$
Since the exponent of $\k$ is not zero (modulo $p$)
, we have found a $c_i'$ which is a primitive $p$'th
root of unity. But we have also proved that we must have 
$p=7$. For if $p\ge13$, there are {\it more than} two elements in 
$G_{k}$, and 
we can make two different choices of $j$, giving conflicting values 
to $c_{k}'$. To sum up, we know that to have collapse we must have 
$p=7$, and some $c_{k}'$ must be a seventh root of unity. Again our 
symmetry argument says that also some $c_{k}''$ must be a seventh root of unity.
The third canonical case is not of interest, since the absolute values are not 
one. We can also easily exclude the second 
canonic case e.g. with the following numerical argument: The imaginary 
part of the seventh power of 
the six values of the components of $c^{(2)}$ are approximately 
$\pm 0.92, \pm 0.94,$ and $\pm 0.41$ rather than 0. So the collapse 
is an internal affair within the first canonical case, which we 
have already studied. This completes 
the proof of the two propositions.
\section{Numerical and asymptotic results}
\label{sec-asymptotics}
\setcounter{equation}{0}

In this section we will study the behavior for large $p$ of the 
solutions $c^{(i)}, i=1,2,3$ defined in Theorem \ref{thm3-1}. We will 
give numerical data leading to educated guesses about this behavior 
(see Remark \ref{observations}
and we will prove quantitative forms of these guesses. 

In Table \ref{sec-asymptotics}.1 below we list the first few primes  $\equiv 1
\pmod{6}$ and corresponding numerical values of 
$A,B,\theta,c_{0}^{(1)},c_{1}^{(1)},$ and $c_{2}^{(1)}$. In Table 
\ref{sec-asymptotics}.2, we give the corresponding information for 
$c^{(2)}$. We will also 
include an indication of the shape of the triangle formed by the 
three complex numbers $c_{0}^{(i)},c_{1}^{(i)},c_{2}^{(i)},\,\, 
(i=1,2),$ reasoning as follows: 

In the 
corresponding situation for simple bi-unimodular sequences of index 
{\it two} (cf. {\cite{pisa}) we have {\it two} complex numbers $c_{0}$ 
and $c_{1}$ on 
the unit circle, and with increasing $p$ their sum tends to zero. A 
natural guess in our situation might therefore be that he sum of
the three numbers tends to zero
 or, equivalently, that the triangle becomes more and more equilateral 
 when $p$ grows. We prefer the latter description. To be able to give
quantitative results we will revive the old noun {\it scalenity}, (cf. 
\cite{oxford}) and give it a precise meaning:
\newtheorem {scal}{Definition}[section] 
\begin{scal} In the complex plane, let $b=(b_{0},b_{1},b_{2})$ be a 
triple of points on a circle $C$
with center $w$. Let $\phi_{i}= \arg(b_{i}-w)$. Let the {\rm 
scalenity} of $b$ be $${\rm scal}(b) = \max_{j} 
\Big|\frac{1}{2}+\cos(\phi_{j+2}-\phi_{j+1})\Big|,$$
(indices counted modulo {\rm 3}).
\end{scal}

\newtheorem{scalformula}[scal]{Remark}
\begin{scalformula}{\rm
\label{scalformula}
Since $\frac{1}{2}=-\cos 2\pi/3$, the triangle with vertices $b$ will be equilateral iff its 
scalenity is zero. Let us now consider the definition of $h_{j}$ (in Proposition 
\ref{prop2-1}). If we take $b=c^{(i)}$ with $i=1$ or 2, we have all $|b_{j}|=1$ and 
thus $w=0$. Hence ${\rm scal}(c^{(1)}) =\frac{1}{2}\max_{j} |1+h_{j}|$, where 
$h_{j}$ is given by (\ref{eq2-3}) with $c$ replaced by $b$.}
\end{scalformula}

\begin{table}[htb]
{\footnotesize
\begin{tabular}{|r|r|r|r|r|r|r|c|} \hline
\multicolumn{8}{|c|}{\normalsize Table \ref{sec-asymptotics}.1 (First canonical case)}\\ \hline
\multicolumn{1}{|c|}{$p$}&\multicolumn{1}{|c|}{$A$}&\multicolumn{1}{|c|}{$B$}&
\multicolumn{1}{|c|}{$\theta$}&\multicolumn{1}{|c|}{$c_{0}^{(1)}$}
&\multicolumn{1}{|c|}{$c_{1}^{(1)}$}&\multicolumn{1}{|c|}{$c_{2}^{(1)}$}
&\multicolumn{1}{|c|}{scal$(c^{(1)})$}\\ \hline
7& 1& 1& 0.4602& $-$0.9010 $-$ 0.4339 i& 0.6235 + 0.7818 i& 
   $-$0.2225 + 0.9749 i& 1.1235\\ \hline 
  13& $-$5& 1& 0.7790& $-$0.4822 $-$ 0.8761 i& 0.3953 + 0.9185 i& 
   $-$0.8132 + 0.5820 i& 0.7132\\ \hline 
  19& 7& 1& 0.2129& $-$0.9528 $-$ 0.3037 i& 0.9838 $-$ 0.1791 i& 
   0.3780 + 0.9258 i& 0.7061\\ \hline 31& 4& 2& 0.4011& $-$0.8023 $-$ 0.5969 i& 
   0.9923 + 0.1235 i& $-$0.0963 + 0.9954 i& 0.5274\\ \hline 
  37& $-$11& 1& 0.9001& $-$0.0604 $-$ 0.9982 i& 0.4630 + 0.8863 i& 
   $-$0.9452 + 0.3265 i& 0.4127\\ \hline 
  43& $-$8& 2& 0.7423& $-$0.3124 $-$ 0.9499 i& 0.7272 + 0.6864 i& 
   $-$0.7742 + 0.6330 i& 0.3792\\ \hline 
  61& 1& 3& 0.5022& $-$0.6466 $-$ 0.7628 i& 0.9759 + 0.2181 i& 
   $-$0.3560 + 0.9345 i& 0.3564\\ \hline 
  67& $-$5& 3& 0.6271& $-$0.4569 $-$ 0.8895 i& 0.8964 + 0.4433 i& 
   $-$0.5999 + 0.8001 i& 0.3170\\ \hline 
  73& 7& 3& 0.3829& $-$0.7843 $-$ 0.6204 i& 0.9988 $-$ 0.0481 i& 
   $-$0.1114 + 0.9938 i& 0.3409\\ \hline 
  79& $-$17& 1& 0.9483& 0.1286 $-$ 0.9917 i& 0.4824 + 0.8759 i& 
   $-$0.9765 + 0.2154 i& 0.3066\\ \hline 
  97& 19& 1& 0.0890& $-$0.9875 $-$ 0.1576 i& 0.7708 $-$ 0.6371 i& 
   0.4823 + 0.8760 i& 0.3137\\ \hline 103& 13& 3& 0.2919& $-$0.8668 $-$ 0.4986 i& 
   0.9629 $-$ 0.2697 i& 0.0653 + 0.9979 i& 0.2937\\ \hline 
  109& $-$2& 4& 0.5556& $-$0.5387 $-$ 0.8425 i& 0.9666 + 0.2561 i& 
   $-$0.4841 + 0.8750 i& 0.2562\\ \hline 
  127& $-$20& 2& 0.8875& 0.0580 $-$ 0.9983 i& 0.6211 + 0.7837 i& 
   $-$0.9411 + 0.3382 i& 0.2464\\ \hline 
  139& $-$23& 1& 0.9731& 0.2257 $-$ 0.9742 i& 0.4899 + 0.8718 i& 
  $-$0.9873 + 0.1588 i& 0.2387\\ \hline 
   151& 19& 3& 0.2290& $-$0.9138 $-$ 0.4062 i& 
    0.9066 $-$ 0.4219 i& 0.1770 + 0.9842 i&0.2452\\ \hline 
  157& $-$14& 4& 0.7212& $-$0.2380 $-$ 0.9713 i& 
    0.8495 + 0.5276 i& $-$0.7655 + 0.6434 i&0.2146\\ \hline 
  163& 25& 1& 0.0683& $-$0.9922 $-$ 0.1248 i& 
    0.7153 $-$ 0.6988 i& 0.4898 + 0.8718 i&0.2411\\ \hline 
  181& 7& 5& 0.4359& $-$0.6932 $-$ 0.7207 i& 
    0.9996 $-$ 0.0270 i& $-$0.2649 + 0.9643 i&0.2092\\ \hline 
\end{tabular}}
\end {table}
}

\noindent
{\footnotesize
\begin{tabular}{|r|r|r|r|r|r|r|c|} \hline
\multicolumn{8}{|c|}{\normalsize Table \ref{sec-asymptotics}.2 (Second canonical case)}\\ \hline
\multicolumn{1}{|c|}{$p$}&\multicolumn{1}{|c|}{$A$}&\multicolumn{1}{|c|}{$B$}&
\multicolumn{1}{|c|}{$\theta$}&\multicolumn{1}{|c|}{$c_{0}^{(2)}$}
&\multicolumn{1}{|c|}{$c_{1}^{(2)}$}&\multicolumn{1}{|c|}{$c_{2}^{(2)}$}
&\multicolumn{1}{|c|}{scal$(c^{(2)})$}\\ \hline
7& 1& 1& 0.4602& 0.8173 + 0.5762 i& $-$0.3890 + 0.9212 i& 0.2804 $-$ 0.9599 i& 
   0.7129\\ \hline 13& $-$5& 1& 0.7790& 0.2469 + 0.9690 i& $-$0.7728 + 0.6346 i& 
   0.6315 $-$ 0.7754 i& 0.9242\\ \hline 
  19& 7& 1& 0.2129& 0.9520 + 0.3061 i& $-$0.4274 + 0.9041 i& 0.0041 $-$ 
  1.0000 i& 
   0.4058\\ \hline 31& 4& 2& 0.4011& 0.8025 + 0.5967 i& $-$0.6855 + 0.7281 i& 
   0.2171 $-$ 0.9761 i& 0.3844\\ \hline 
  37& $-$11& 1& 0.9001& $-$0.2907 + 0.9568 i& $-$0.9939 $-$ 0.1106 i& 
   0.8980 $-$ 0.4399 i& 0.6830\\ \hline 43& $-$8& 2& 0.7423& 0.1847 + 0.9828 i& 
   $-$0.9804 + 0.1971 i& 0.7022 $-$ 0.7120 i& 0.5126\\ \hline 
  61& 1& 3& 0.5022& 0.6436 + 0.7653 i& $-$0.8671 + 0.4982 i& 
   0.3730 $-$ 0.9278 i& 0.3232\\ \hline 67& $-$5& 3& 0.6271& 0.4178 + 0.9085 i& 
   $-$0.9589 + 0.2837 i& 0.5632 $-$ 0.8263 i& 0.3571\\ \hline 
  73& 7& 3& 0.3829& 0.7932 + 0.6089 i& $-$0.7822 + 0.6231 i& 
   0.1987 $-$ 0.9801 i& 0.2661\\ \hline 
  79& $-$17& 1& 0.9483& $-$0.4295 + 0.9031 i& $-$0.8761 $-$ 0.4821 i& 
   0.9657 $-$ 0.2596 i& 0.4408\\ \hline 
  97& 19& 1& 0.0890& 0.9887 + 0.1498 i& $-$0.4791 + 0.8778 i& 
   $-$0.2383 $-$ 0.9712 i& 0.2383\\ \hline 
  103& 13& 3& 0.2919& 0.8757 + 0.4828 i& $-$0.7239 + 0.6899 i& 
   0.0609 $-$ 0.9981 i& 0.2327\\ \hline 
  109& $-$2& 4& 0.5556& 0.5288 + 0.8487 i& $-$0.9508 + 0.3097 i& 
   0.4751 $-$ 0.8799 i& 0.2600\\ \hline 127& $-$20& 2& 0.8875& $-$0.2473 + 0.9689 i& 
   $-$0.9083 $-$ 0.4183 i& 0.9253 $-$ 0.3794 i& 0.3193\\ \hline 
  139& $-$23& 1& 0.9731& $-$0.4663 + 0.8846 i& $-$0.7821 $-$ 0.6231 i& 
   0.9836 $-$ 0.1804 i& 0.3135\\ \hline
151& 19& 3& 0.2290& 0.9203 + 0.3913 i& 
    $-$0.6848 + 0.7287 i& $-$0.0415 $-$ 0.9991 i&0.1997\\ \hline 
  157& $-$14& 4& 0.7212& 0.1708 + 0.9853 i& 
    $-$0.9969 $-$ 0.0789 i& 0.7371 $-$ 0.6757 i&0.2520\\ \hline 
  163& 25& 1& 0.0683& 0.9929 + 0.1190 i& 
    $-$0.4865 + 0.8737 i& $-$0.2912 $-$ 0.9567 i&0.1941\\ \hline 
  181& 7& 5& 0.4359& 0.7010 + 0.7132 i& $-$0.9002 + 0.4355 i& 
    0.2959 $-$ 0.9552 i&0.1824\\ \hline

\end{tabular}}
\bigskip

We present the corresponding values for $c^{(3)}$ in Table 
\ref{sec-asymptotics}.3. Since these values are real,
we will save some space and we use this for giving the information 
also in another form, namely $ \frac{c_{j}^{(3)}}{\sqrt{p}}$, 
which should shed some light on the surprising behaviour of the 
components.
\vskip5mm \noindent
{\footnotesize
\begin{tabular}{|r|r|r|r|r|r|r|r|r|c|} \hline
\multicolumn{10}{|c|}{\normalsize Table \ref{sec-asymptotics}.3 (Third canonical 
case)}
\\ \hline
\multicolumn{1}{|c|}{$p$}&\multicolumn{1}{|c|}{$A$}&\multicolumn{1}{|c|}{$B$}&
\multicolumn{1}{|c|}{$\theta$}&\multicolumn{1}{|c|}{$c_{0}^{(3)}$}
&\multicolumn{1}{|c|}{$c_{1}^{(3)}$}&\multicolumn{1}{|c|}{$c_{2}^{(3)}$}
&\multicolumn{1}{|c|}{$c_{0}^{(3)}/\sqrt{p}$}
&\multicolumn{1}{|c|}{$c_{1}^{(3)}/\sqrt{p}$}&\multicolumn{1}{|c|}{$c_{2}^{(3)}/\sqrt{p}$}
\\ \hline

7&1&1& 0.4602& $-$1.2221& 9.4127& 2.7389&  $-$0.4619& 3.5577& 1.0352\\ \hline 
  13& $-$5& 1& 0.7790& $-$1.4201& $-$14.6415& 2.1601& $-$0.3939& 
  $-$4.0608& 0.5991\\ \hline
  19& 7& 1& 0.2129& $-$2.2521& 8.4655& 4.8488& $-$0.5167& 1.9421& 1.112\\ \hline
  31& 4& 2& 0.4011& $-$2.8168& 17.2938& 4.6888&$-$0.5059& 3.1061& 0.8421\\ \hline
  37& $-$11& 1& 0.9001& $-$3.0328& $-$7.1015& 2.8445&  $-$0.4986& $-$1.1675& 0.4676\\ \hline
  43& $-$8& 2& 0.7423& $-$3.2558& $-$20.3776& 3.6527& $-$0.4965& $-$3.1075& 0.557\\ \hline
  61& 1& 3& 0.5022& $-$4.0014& 50.9574& 5.4586& $-$0.5123& 6.5244& 0.6989\\ \hline
  67& $-$5& 3& 0.6271& $-$4.2289& $-$95.9688& 5.0005&$-$0.5166&$-$11.7245& 0.6109\\ \hline
  73& 7& 3& 0.3829& -4.4100& 25.6091& 6.6407& -0.5162& 2.9973& 0.7772\\ \hline 
  79& -17& 1& 0.9483& -5.6126& -8.9422& 4.1623& -0.6315& -1.0061& 0.4683\\ \hline 
  97& 19& 1& 0.0890& -5.0982& 13.5365& 10.4124& -0.5176& 1.3744& 1.0572\\ \hline 
  103& 13& 3& 0.2919& -5.2556& 21.9500& 8.4142& -0.5179& 2.1628& 0.8291\\ \hline 
  109& -2& 4& 0.5556& -5.5068& 337.8101& 6.6180& -0.5275& 32.3563& 0.6339\\ \hline 
  127& -20& 2& 0.8875& -7.1148& -14.4873& 5.5161& -0.6313& -1.2855& 0.4895\\ \hline 
  139& -23& 1& 0.9731& -8.4417& -11.5986& 5.6060& -0.7160& -0.9838& 0.4755\\ \hline 
  151& 19& 3& 0.2209& -6.3543& 22.1314& 10.5843& -0.5171& 1.8010& 0.8613\\ \hline 
  157& -14& 4& 0.7212& -7.1235& -36.5459& 6.8056& -0.5685& -2.9167& 0.5431\\ \hline 
  163& 25& 1& 0.0683& -6.5753& 16.4057& 13.3294& -0.5150& 1.2850& 1.0440\\ \hline 
  181& 7& 5& 0.4359& -7.0841& 58.2887& 9.2448& -0.5266& 4.3326& 0.6872\\ \hline 

\end{tabular}}

\bigskip
Our observations are summarized in the following remark:
\newtheorem{observations}[scal]{Remark}
\begin{observations}{\rm Our numerical observations and our results are of five kinds:

\noindent (1)  For each large $p$, the first and second canonical solutions 
are approximately symmetric to each other w.r.t. the origin. 

\noindent (2) Even though two large primes may 
be close to each other without their canonical solutions being close, 
large primes with approximately the same 
$\theta$ will have approximately the same first canonical solutions 
and approximately the same second canonical solutions (even if 
the primes are not close to each other). 

\noindent (3) For large $p$, the first and second canonical solution each forms an approximately 
equilateral triangle. 

\noindent(4) For large $p$, the approximate {\it positions} 
of the nearly equilateral triangles are simple functions of $\theta$.

\noindent (5) If $p$ is large, then all components of $|c^{(3)}|$ are 
large. If in addition $|A|$ is small, that is if $\theta$ is close to $\pi/6$, then 
$|c_{1}^{(3)}|$ is very large.
\label{observations}
}\end{observations}
To make it easier to guess {\it quantitative} results (making 
``approximately'' more precise in Remark \ref{observations}) we 
present a few more numerical results in Table \ref{sec-asymptotics}.4.

\begin{table}[htb]  
{\footnotesize
\begin{tabular}{|r|r|r|r|r|r|r|r|c|} \hline
\multicolumn{9}{|c|}{\normalsize Table \ref{sec-asymptotics}.4 (Large
primes, 
close in size vs. close in $\theta$-value)}\\ \hline
\multicolumn{1}{|c|}{$p$}&\multicolumn{1}{|c|}{$A$}&\multicolumn{1}{|c|}{$B$}&
\multicolumn{1}{|c|}{$\theta$}&\multicolumn{1}{|c|}{$\arg(c_{0}^{(1)})$}
&\multicolumn{1}{|c|}{$2\theta-\pi$}
&\multicolumn{1}{|c|}{$\arg(c_{0}^{(2)})$}&\multicolumn{1}{|c|}{$2\theta$}
&\multicolumn{1}{|c|}{scal($c^{(1)})$}\\ \hline
1003273&973&337& 0.354542&$-$2.43320&$-$2.43251&0.70803&0.709084&0.002810\\ \hline  
1003279&1993&39&0.033775&$-$3.07411&$-$3.07404&0.06742&0.067555&0.002995\\ \hline 
100205473&9733&3367&  0.354372& $-$2.43292&$-$2.43285& 0.70864& 
0.708744&0.000281\\ \hline
  \end{tabular}}
  \end{table}
 
From Table \ref{sec-asymptotics}.4 it seems that
 ``approximately''means agreement in approximately $\frac{n}{2}$ 
decimals. Thus quantitative results in terms of 
$O(\frac{1}{\sqrt{p}})$ might seem 
plausible. In our quantitative results we will use the maximum norm 
to measure distances in $\bC^{\,3}$. We will also need a name for the equilateral 
``limit'' triangle hinted atin Remark \ref{observations} (4), 
hopefully visible in Tables \ref{sec-asymptotics}.2 and 
\ref{sec-asymptotics}.3, and present in columns 5 and 8 of Table \ref{sec-asymptotics}.4. 
Thus we make the following two definitions:
\newtheorem {norm}[scal]{Definition} 

\begin{norm}
\label{norm} Let $a=(a_{0},a_{1},a_{2})\in\bC^{\,3}$, then we define 
$\|a\|=\max(|a_{0}|,|a_{1}|,|a_{2}|).$
\end{norm}
\newtheorem {limit}[scal]{Definition}
\begin{limit}
\label{limit}Let $p$ be a prime $\equiv 1 \pmod{6}$ and let 
$\theta=\frac13 \Arccos\left(\frac{A}{2\sqrt{p}}\right)$, where 
$4p=A^2+27B^2$ and $A\equiv 1 \pmod{3}$.
We denote by $d=d(p)=(d_{0},d_{1},d_{2})$ the (equilateral) triangle for 
which 
$$d_{j}=\exp\Bigl(2i(\theta-\frac{2j\pi}{3})\Bigr),\,\,\,j=0,1,2.$$
\end{limit}

We will now state four quantitative results for the first and second 
canonical cases, where Proposition 
$\ref{sec-asymptotics}.j$  for $j = 6,\ldots,9$ is of the kind $(j-5)$ listed 
in Remark \ref{observations}. (The discusion of kind (5) starts after 
Corollary \ref{cor} below.)
\newtheorem  {c1+c2}[scal]{Proposition}
\begin{c1+c2} Let $p$ be a prime $\equiv 1 \pmod{6}$, and let  
$c^{(1)}$ and $c^{(2)}$ be the corresponding first and second canonical solution.
Then 
$$ \|c^{(1)}+c^{(2)}\|\le \frac{36}{5\sqrt{p}}.$$ 
\label{c1+c2}
\end{c1+c2}

\newtheorem  {p1p2}[scal]{Proposition}
\begin{p1p2} Let $p'$and $p''$ be primes $\equiv 1 \pmod{6}$, let  
$\theta'$ and $\theta''$ be their respective $\theta$-values and let 
$c'$ and $c''$ be their respective first canonical solutions. Then
$$ \|c'-c''\|\le 2|\theta'-\theta''|+\frac{3}{\sqrt{p'}}+ 
\frac{3}{\sqrt{p''}}.$$ 
The same result, with the constants {\rm 3} replaced by 
$\frac{21}{5}$, holds if $c'$ and $c''$ are the respective {\rm second} canonical solutions.
\label{p1p2}
\end{p1p2}

\newtheorem  {smallscal}[scal]{Proposition}
\begin{smallscal} Let $p$ be a prime $\equiv 1 \pmod{6}$, and let  
$c^{(1)}$ and $c^{(2)}$ be the corresponding first and second canonical solution. Then 
$$ {\rm scal}(c^{(1)})\le \frac{7}{2\sqrt{p}} \,\,\,\,and\,\,\,\, 
{\rm scal}(c^{(2)})\le \frac{21}{5\sqrt{p}}.$$ 
\label{smallscal}
\vskip-8mm
\end{smallscal}
\newtheorem  {limc1c2}[scal]{Proposition}
\begin{limc1c2}Let $p$ be a prime $\equiv 1 \pmod{6}$, let  
$c^{(1)}$ and $c^{(2)}$ be the corresponding first and second canonical 
solution, and let $d$ be as in Definition {\rm\ref{limit}}. Then
$$ \|c^{(1)}+d\|\le \frac{3}{\sqrt{p}}\,\,{\rm and}\,\,
\|c^{(2)}-d\|\le \frac{21}{5\sqrt{p}}.$$
\label{limc1c2}
\vskip-8mm
\end{limc1c2}
\newtheorem{constants}[scal]{Remark}
\begin{constants}{\rm
The constants in these propositions are not best possible but are 
chosen as compromises to make the proofs less cumbersome. Even if we 
restrict our claims to hold only for $p>M$ for some large $M$, the 
constants cannot always be significantly improved. For instance,
for $p=10^{10}+279$ we have $\|c^{(2)}-d\|\approx 4/\sqrt{p}.$
For a kind of ``best possible'', result, see Remark \ref{best}.
\label{constants}
}\end{constants}

Since to each number $\theta$ (in the interval $[0,\pi/3]$) there 
corresponds
at most one $p$, it does not make sense to consider a sequence of $p$:s 
with a common $\theta$. However, Proposition \ref{limc1c2} 
obviously has the following corollary, where we have used the notation  
$\theta(p)$, $c^{(2)}(p)$ and $c^{(1)}(p)$ for the values of $\theta$ and the
first and second canonical solutions corresponding to $p$:
\newtheorem  {cor}[scal]{Corollary}
\begin{cor}
\label{cor}
Let $\theta_{0}$ be a real number in the interval $[0,\pi/3]$.
Denote by $d=(d_{0},d_{1},d_{2})$ the (equilateral) triangle for 
which 
$\displaystyle 
d_{j}=\exp\Bigl(2i(\theta-\frac{2j\pi}{3})\Bigr),\,\,\,j=0,1,2.$
Let $\{p_{n}\}_{1}^\infty$ be a sequence of primes $\equiv 1 \pmod{6}$   
going to infinity in such a way that $\lim_{n\rightarrow \infty} 
\theta(p_{n})=\theta_{0}$. Then $\lim_{n\rightarrow 
\infty}c^{(1)}(p_{n})=-d$ and $\lim_{n\rightarrow 
\infty}c^{(2)}(p_{n})=d$.
\end{cor}

Before proving our four Propositions we will comment item (5) of 
Remark \ref{observations}. In Table \ref{sec-asymptotics}.5 we present some more 
numerical values with focus on $\theta$-values close to $0, \pi/6$ 
and $\pi/3$.
\vskip5mm\noindent
{\footnotesize
\begin{tabular}{|r|r|r|r|r|r|r|} \hline
\multicolumn{7}{|c|}{\normalsize Table \ref{sec-asymptotics}.5 (Third 
canonical case for large primes)}\\ \hline
\multicolumn{1}{|c|}{$p$}&\multicolumn{1}{|c|}{$A$}&\multicolumn{1}{|c|}{$B$}&
\multicolumn{1}{|c|}{$\theta$}&\multicolumn{1}{|c|}{$c_{0}^{(3)}/\sqrt{p}$}
&\multicolumn{1}{|c|}{$c_{1}^{(3)}/\sqrt{p}$}&\multicolumn{1}{|c|}{$c_{2}^{(3)}/\sqrt{p}$}
\\ \hline
67 521 601 729&$-$2&100 016& 0.523600&$-$0.577349&779 550.5&0.577353\\ \hline  
67 544 557 351&1&100 033&0.523598&$-$0.577347&194 920.0&0.577353\\ \hline 
250 004 500 027&1 000 009&1&0.000002&$-$0.500000&1.000007&1.000001 \\ \hline
250 018 500 349&$-$1 000 037&1&1.047196&$-$0.999995&$-$0.999997&0.500000 \\ \hline
  \end{tabular}}\break
\vskip5mm
 These and other numerical results make it plausible that ``large'' in item (5) 
 of Remark 
 \ref{observations}, may be specified to mean ``not much smaller 
 than $\frac{1}{2}\sqrt{p}$ '', but it seems
difficult to find $\theta$-independent estimates of ``convergence 
rate'' for the third canonical case. We are now ready to state a 
proposition:

\newtheorem{third}[scal]{Proposition}
\begin{third} If 
 $\{p_{n}\}_{1}^\infty$ is any sequence of primes $\equiv 1 \pmod{6}$   
going to infinity, then (with obviuos notation) for $i=0,1,$ and {\rm 
2},
\begin{equation}
\label{liminf} 
\liminf_{n\rightarrow 
\infty}
\frac{|c_{i}^{(3)}(p_{n})|}{\sqrt{p_{n}}} \ge 0.5. 
\end{equation}
\label{third}
\end{third}

We remark that this proposition implies that for every normalized 
$x=(1,x_{1}\ldots,x_{p-1})\in\bR^p$  of index 
3 coming from the third canonical case for
a large $p$, either all $|x_{j}|,\,\, j\ne 0,$ are large or they are all small
(leaving the canonical case via the transformations mentioned in 
Theorem \ref{thm3-1} and leaving the simple case via Definition 
\ref{general}). 

We will now prove our five propositions. 
\medskip 

\noindent{\bf Proof of Propositions \ref{c1+c2} and \ref{c1+c2}}
Proposition \ref{c1+c2} follows from Proposition \ref{limc1c2}
via a straightforward application of the triangle inequality. 

Similarly, Proposition \ref{p1p2} follows from Proposition \ref{limc1c2} via 
the triangle inequality and the inequality
$|\exp(2 i\phi') - \exp(2 i\phi'')|\le 2|\phi'-\phi''|.$ 
\medskip

\noindent{\bf Proof of Proposition \ref{smallscal}}. From \eqref {eq2-3} 
and \eqref {eq2-5} and Remark \ref{scalformula} we get 
\begin{equation}
\label{scal1}
\sqrt{p}\,{\rm scal}\,\,(c^{(1)})=\tfrac12\sqrt{p}\max_{j}|h_{j}
+1|\le\tfrac12\sqrt{p}\,\Big(|1+\xi_{1}^{(1)}|\,+\,|\eta_{1}^{(1)}|\Big)=
\frac{ 3\sqrt{p}(p-A) + p(6p-24)}{2(p^2-3p-A)}.
\end{equation}
Since the right-hand side of of \eqref{scal1} is a decreasing function of $A$ and  
$A>-2\sqrt{p}$, we get 
\begin{equation}
\label{scal1no2}
\sqrt{p}\,{\rm scal}\,\,(c^{(1)})< \frac{ 3\sqrt{p}(p+2\sqrt{p}) + 
p(6p-24)}{2(p^2-3p+2\sqrt{p})}= 3 + 
\frac{3(\sqrt{p}-2)}{2(\sqrt{p}-1)^2}.
\end{equation}
The last member 
of \eqref{scal1no2} as a function of $p$ is decreasing for $p>9$ and 
takes values $< 3.4$ for $p=7$ and 13. This completes the proof of 
the first part of the proposition.

For the second part, we will use \eqref{eq2-6} and the identity 
\[
\Big(\sqrt{p}\sqrt{p+4A+16}+p+2\Big)\Big(\sqrt{p}\sqrt{p+4A+16}-p-2\Big)=
4(Ap+3p-1)
\]
to give the counterpart of \eqref{scal1} the form
\[
\sqrt{p}\,{\rm scal}\,\,(c^{(2)})\le\tfrac12\sqrt{p}\,\Big(|1+\xi_{2}^{(1)}|\,+\,|\eta_{2}^{(1)}|\Big)=
\frac{ 3(2p+\sqrt{p})}{p+2+\sqrt{p}\sqrt{p+4A+16}}.
\]
We can again take $A=-2\sqrt{p}$. The resulting 
expression is easily seen to be $<4$ for $p>100$ and for the 
remaining $p$ we enter the true value of $A$ (given 
in Table \ref{sec-asymptotics}.1) to get 
a maximum  $\approx 4.1966$ for $p=37$. This completes the proof.

\newtheorem{best}[scal]{Remark}
\begin{best}{\rm From \eqref {scal1no2} we easily get the following 
result: For each $\epsilon$ with $0<\epsilon<1$ we have
\[
{\rm scal}(c^{(1)})\le \frac{3+\epsilon}{\sqrt{p}} \,\,\,\,if\,\,\,\, 
p>\Big(2+\frac{3}{2\epsilon}\Big)^2,
\]
which could be contrasted with the fact that for $p=$ 10 002 900 217 we 
have $\sqrt{p}\,\, {\rm scal}(c^{(1)})\approx 3.000015.$ 
\label{best}
}\end{best}

In the proof of Proposition 
\ref{limc1c2} we will work with $\alpha$, $\beta$ and $\gamma$ as 
given in Theorem \ref{thm3-1} and $\rho$, $\sigma$, and $\tau$ as given in 
Lemma \ref{lemma3-2}. We will use the following lemma:
\newtheorem  {alphabetagamma}[scal]{Lemma}
\begin{alphabetagamma} Let $b'=(b_{0}',b_{1}',b_{2}')\in\bC^{\,3}$
 and $b''=(b_{0}'',b_{1}'',b_{2}'')\in\bC^{\,3}$ be given by
 \[
b_j'=\rho'+\sigma'\cos\left(\theta-\frac{2\pi}{3}j\right)+\tau'
\sin\left(\theta-\frac{2\pi}{3} j\right),\quad j=0,1,2,
\]
\begin{equation}
\label{bbis}
b_j''=\rho''+\sigma''\cos\left(\theta-\frac{2\pi}{3}j\right)+\tau''
\sin\left(\theta-\frac{2\pi}{3} j\right),\quad j=0,1,2,
\end{equation}
where $\theta\in\bR$ and 
$\rho',\sigma',\tau',\rho'',\sigma'',\tau''\in \bC$,
Then 
\[\|b'-b''\|\le 
|\rho'-\rho''|+\sqrt{|\sigma'-\sigma''|^2+|\tau'-\tau''|^2}.
\]
\label{alphabetagamma}
\end{alphabetagamma}
The proof of Lemma \ref{alphabetagamma} is a straightforward 
application of the triangle inequality, the Cauchy inequality, and the 
identity $\cos^2+\sin^2=1$.

\medskip 

\noindent{\bf Proof of Proposition \ref{limc1c2}}  In Lemma 
\ref{alphabetagamma} we take $b'=c^{(1)}$ and $b''=-d$ (cf. Definitions 
\ref{canonic} and \ref{limit}). Then 
$\rho'=\alpha^{(1)},\sigma'=\beta^{(1)},\tau'=\gamma^{(1)}$ as given 
in \eqref{eq3-3}, whereas $\rho''=0,\sigma''=-\cos 3\theta-i 
\sin 3\theta$, and $\tau''=-\sin 3\theta+ i \cos 3\theta$, as is 
easily checked by introducing these values in \eqref{bbis}
and applying the addition theorems for sine and cosine. Since $\displaystyle\cos 
3\theta=\frac{A}{2\sqrt{p}}$ and $\displaystyle\sin 
3\theta=\frac{3B\sqrt{3}}{2\sqrt{p}}$, Lemma \ref{alphabetagamma} 
shows that for the proof of the first half of Proposition 
\ref{limc1c2} it only remains to check that with $\alpha^{(1)}, 
\beta^{(1)},$ and $\gamma^{(1]}$ as in \eqref{eq3-3} we have
\begin{equation}
\sqrt{p}\,|\alpha^{(1)}|+\sqrt{p\Big|\beta^{(1)}+
\frac{A+ i3 \sqrt{3}B}{2\sqrt{p}}\Big|^2+p\Big|\gamma^{(1)}+
\frac{3\sqrt{3}B-iA} {2\sqrt{p}}\Big|^2}\le 3. 
\label{c1vsd}
\end{equation}
Introducing the values of $\alpha^{(1)}, \beta^{(1)},$ and 
$\gamma^{(1]}$ and replacing $3B\sqrt{3}$ by $\sqrt{4p-A^2}$ 
we can after some calculation treat the first term of the left 
member of \eqref{c1vsd} as follows
\begin{equation}
\sqrt{p}\,|\alpha^{(1)}|=\sqrt{\frac{p^2-A 
p}{p^2-3p-A}}<\sqrt{\frac{p^2+2p\sqrt{p}}
{p^2-3p+2\sqrt{p}}}=\sqrt{\frac{p(2+\sqrt{p})}{(p-3)\sqrt{p}+2}},
\label{firstterm}
\end{equation}
where the estimate comes from the facts that the second term of 
\eqref{firstterm} is a decreasing function of $A$ and that 
$A>-2\sqrt{p}$. Let us denote by $Q$ the expression under the big 
root sign in \eqref{c1vsd}. Since the last member of \eqref{firstterm} is a 
decreasing function of $p$ with a value $< 1.5$ for $p=31$, we can prove
\eqref{c1vsd} for $p\ge 31$ 
by checking that 
\begin{equation}
Q\le (3-1.5)^2=2.25\,\, {\rm for}\,\, p\ge 31
\label{Qsmall}
\end{equation}
 Treating Q in the same way as we did with first term of the left 
member of \eqref{c1vsd} we find
\begin{equation}
Q=\frac{2p^3-(A+6)p^2+2A p-(2p-A-4)\sqrt{p^4-4p^3}}{p^2-3p-A}.
\label{Q}
\end{equation}
Using a Taylor formula with rest term we have 
$$\sqrt{p^4-4p^3}=p^2\Big(1-\frac{4}{p}\Big)^{\frac{1}{2}}=p^2-2p-2+R_{3},$$
where $-\frac{6}{p}<R_{3}<0$ (since $p>31$).
Introducing this in \eqref{Q} we get
\begin{equation}
Q=\frac{2p^2-4p-2A-8-(2p-A-4)R_{3}}{p^2-3p-A}<\frac{2(p^3-2p^2+Ap-2p-3A-12)}{p(p^2-3p-A)}.
\label{QTaylor}
\end{equation}
Since the last member of \eqref{QTaylor} is an increasing function 
of $A$ we can estimate it with its value for $A=2\sqrt{p}$, which is 
a decreasing function of $p$ and thus not larger than its value for 
$p=31$, 
which turns out to be $\approx 2.07$ in agreement with \eqref{Qsmall}. Finally, we 
check numerically the value of $\sqrt{p} \|c^{(1)}+d\|$ for $p=7,13,$ and $19$. 
We find 2.59, 2.31, and 1.91, which are all $<3$. 
This completes the proof of the first half of the proposition. 

For the second part of the proof we proceed in the same way but let MATHEMATICA help us 
to get a good start, namely by telling us that defining 
$m(p)=\sqrt{p}\, \|c^{(2)}-d\|$ we have $m(p)\le m(43)< 4.1$ if $p < 
10 000$. We get after some calculation
\[
\sqrt{p}\,|\alpha^{(2)}|=\sqrt{\frac{2p}{2+p+\sqrt{p}\sqrt{p+4A+16}}}\le
\sqrt{\frac{2p}{2+p+\sqrt{p}\sqrt{p-8\sqrt{p}+16}}}<
\frac{100}{99}
\]
if $p>10 000$. Thus to complete the proof is enough to prove that
\begin{equation}
\label {less10.17}
p\Big|\beta^{(2)}-
\frac{A+ i3 \sqrt{3}B}{2\sqrt{p}}\Big|^2+p\Big|\gamma^{(2)}-
\frac{3\sqrt{3}B-iA} {2\sqrt{p}}\Big|^2\le 10.17
\end{equation}
 if $p>10 000$, e.g. by 
 proving that the first term of \eqref{less10.17} is $< 1.07$ and the 
 second term is $< 9.1.$. This can be done as in the proof of the first 
 part, using \eqref{eq3-4}. Just as we have studied functions 
 of $A$ restricted to the interval $|A|<2\sqrt{p}$, we will now 
 with the help of \eqref{Atouv} and \eqref{Btouv} 
write the left member of \eqref{less10.17} as a function of $u$ 
and $v$, where $|u-4|<v<u+4$. Again a certain square root
can be estimated with a Taylor formula. We leave the details to the reader.
}}
\medskip

\noindent{\bf Proof of Proposition \ref{third}}
Inspired by the first two rows of Table \ref{sec-asymptotics}.5 we 
expect infinities near $\theta=\pi/6$, and thus, to avoid zeros in 
the denominator, we ``turn everything upside down''. Thus we want to 
prove that 
$$\limsup_{n\rightarrow \infty}
\frac{\sqrt{p_{n}}}{|c_{i}^{(3)}(p_{n})|} \le 2.$$
Suppose this is not true. Then (by taking subsewuences if 
needed) we can find a sequence $\{p_{n}\}_{1}^\infty$ 
of primes $\equiv 1 \pmod{6}$ going to infinity, 
such that
\begin{equation}
\label{lim} 
\lim_{n\rightarrow 
\infty}
\frac{\sqrt{p_{n}}}{c_{j}^{(3)}(p_{n})}{\sqrt{p_{n}}} =l_{j}, 
\end{equation}
where these limits 
exist (finite or $+\infty$) and $|l_{j}|>2$ for at least 
one $j$ (0,1, or 2). Since the interval $[0,\pi/3]$ is compact, we can 
by again taking a subsequence (keeping the notation 
$\{p_{n}\}_{1}^\infty$) arrange that $\theta_{0}=\lim_{n\rightarrow 
\infty}\theta(p_{n})$ exists. 
Starting from \eqref{eq3-5} we replace $A$ by 
$2\sqrt{p}\,\cos3\theta$ and $B$ by $2\sqrt{p}\,\sin{3\theta}/\sqrt{27}$. 
Introducing the resulting expressions for $\alpha^{(3)}, \beta^{(3)},$ and 
$\gamma^{(3)}$ in \eqref{eq3-2}, we get $\sqrt{p}/c^{(3)}$ as a function of 
$p$ and $\theta$, which we denote by 
$q(p,\theta)$. We now fix
$\theta=\theta_{0}$ and study $q(p,\theta_{0})$ as a function of $p$ 
when $p\rightarrow \infty$. Estimating various square roots with a 
Taylor formula, we get after a considerable amount of calculation:
%
\[
\label{limq} 
\lim_{p\rightarrow 
\infty}q(p,\theta_{0})=\Big(-2\cos\theta_{0}\,\,,
\,\,-2\sin(\theta_{0}-\pi/6)\,\,,\,\,
2\sin(\theta_{0}+\pi/6)\Big).
\]
A simple continuity argument (w.r.t. $\theta(p_{n}$ and
$\theta_{0}$) shows that with $l_{j}$ from \eqref{lim} we have
\begin{equation} 
l_{0}=-\cos\theta_{0}\,\,,
\,\,l_{1}=-2\sin(\theta_{0}-\pi/6)\,\,,\,\,
l_{2}=2\sin(\theta_{0}+\pi/6).
\label{l}
\end{equation}

This is a contradiction, since we have supposed that $|l_{j}|>2$ 
for at least one $j$. We have thus completed the proof and also 
substantiated the ``very large'' part of item (5) of 
Remark \ref{observations} (take $l_{1}$ from \eqref{l} and consider 
$|1/l_{1}|$ for $\theta_{0}$ close to $\pi/6$).





\end{document}